\documentclass[11pt]{article}

\usepackage[a4paper,margin=1in]{geometry}
\usepackage[T1]{fontenc}
\usepackage{lmodern}
\usepackage{amsmath,amssymb,amsthm,mathtools}
\usepackage{microtype}
\usepackage{enumitem}
\usepackage{booktabs}
\usepackage{xcolor}
\usepackage[colorlinks=true,linkcolor=blue!55!black,citecolor=blue!55!black,urlcolor=blue!55!black]{hyperref}
\usepackage[nameinlink,capitalize,noabbrev]{cleveref}
\usepackage{fancyhdr}
\newtheorem{theorem}{Theorem}[section]
\newtheorem{proposition}[theorem]{Proposition}
\newtheorem{lemma}[theorem]{Lemma}
\newtheorem{corollary}[theorem]{Corollary}
\newtheorem{remark}[theorem]{Remark}

\crefname{assumption}{Assumption}{Assumptions}
\Crefname{assumption}{Assumption}{Assumptions}
\crefname{theorem}{Theorem}{Theorems}
\Crefname{theorem}{Theorem}{Theorems}
\crefname{proposition}{Proposition}{Propositions}
\Crefname{proposition}{Proposition}{Propositions}
\crefname{lemma}{Lemma}{Lemmas}
\Crefname{lemma}{Lemma}{Lemmas}
\crefname{corollary}{Corollary}{Corollaries}
\Crefname{corollary}{Corollary}{Corollaries}
\crefname{remark}{Remark}{Remarks}
\Crefname{remark}{Remark}{Remarks}

\newcommand{\C}{\mathbb C}
\newcommand{\R}{\mathbb R}
\newcommand{\N}{\mathbb N}
\newcommand{\Pp}{\mathbb P}
\newcommand{\Ee}{\mathbb E}

\newcommand{\dist}{\operatorname{dist}}

\newcommand{\rank}{\operatorname{rank}}
\newcommand{\Span}{\operatorname{span}}

\newcommand{\diag}{\operatorname{diag}}

\newcommand{\norm}[1]{\left\lVert #1\right\rVert}

\definecolor{darkblue}{rgb}{.2, 0.2,.8}
\definecolor{darkgreen}{rgb}{0,0.5,0.3}
\definecolor{darkred}{rgb}{.8, .1,.1}

\newcommand{\nto}{n\to\infty}
\renewcommand{\P }{{\mathbb P}}

\title{Spectrum of non-Hermitian self-normalized random matrices with heavy tails}

\author{Johannes Heiny \and Michail Louvaris}

\date{\today}

\begin{document}

\maketitle

\begin{center}
\small
Department of Mathematics, KTH Royal Institute of Technology, Stockholm, Sweden\\
\texttt{heiny@kth.se}  \qquad \texttt{louvaris@kth.se}
\end{center}

\maketitle
\begin{abstract}
We study non-Hermitian random matrices $Y_n$ obtained by normalizing each
row of an $n\times n$ random matrix with i.i.d.\ entries to have Euclidean norm
one. The modulus of the entries is assumed to be in the domain of attraction of an $\alpha$-stable distribution with $\alpha\in[0,2]$, where $\alpha=0$ is interpreted as slowly varying tails.
%have a regularly varying radial tail $\mathbb P(|X|>t)=t^{-\alpha}L(t)$ with index $ \alpha\in[0,2]$, where $L$ is a slowly varying function. 
We determine the limiting singular value
and eigenvalue distributions of $Y_n$.

For $\alpha\in(0,2)$, we prove convergence of the empirical
singular value measures of every shifted matrix $Y_n-zI$ to
deterministic laws characterized by the root resolvent of a
column-normalized, directed Poisson Weighted Infinite Tree (PWIT), see \cite{BCC}.
Under a bounded density assumption, we establish the convergence of the empirical spectral distribution of $Y_n$ to a deterministic measure $\mu_\alpha^{\mathrm{SN}}$ on $\C$. $\mu_\alpha^{\mathrm{SN}}$ is rotationally invariant, is supported in
the closed unit disk with outer support radius one, has no point atoms,
and assigns no mass to the unit circle. The laws \(\mu_\alpha^{\mathrm{SN}}\) interpolate between \(\delta_0\) at \(\alpha=0\) and the circular law at \(\alpha=2\). These also arise as the limiting eigenvalue distributions of \(Y_n\) for entries with slowly varying tails, $\alpha=0$, and for centered entries with bounded density in the normal domain of attraction, respectively.

As an application of our results and based on the upper bound provided by \cite{franccois2026characteristic}, we derive the convergence of the spectral radius of $Y_n$ to $1$ in probability for $\alpha \in (0,2]$. 
Our approach yields a general criterion for the outer support radius
of rotationally invariant measures characterized by a fresh cavity
recursion. Under explicit trace and radial-recovery assumptions, this
radius is determined by the fractional moments of the recursion weights.
\end{abstract}
\section{Introduction}

The celebrated circular law describes the limiting eigenvalue distribution of suitably normalized random matrices with independent entries of finite variance; see, for instance, \cite{TaoVu2010,bordenave2012around}. In the infinite variance regime, deterministic normalization leads instead to a family of non-Gaussian spectral laws, whose analysis for non-Hermitian heavy-tailed matrices was developed by Bordenave, Caputo and Chafa\"i \cite{BCC}. A natural alternative is to normalize each row by its own Euclidean norm. This fixes the length of every row, but introduces  dependence between the entries within a row. While this dependence is very weak in light-tailed scenarios, it becomes significant under heavy tails. The purpose of this paper is to describe the limiting spectral behavior of such self-normalized matrices throughout the heavy-tailed regime.

Self-normalization arises naturally in multivariate statistics, most notably through sample correlation matrices. If a data matrix is normalized row by row, then the associated Gram matrix is a sample correlation matrix. In contrast with sample covariance matrices, sample correlation matrices are invariant under coordinatewise rescaling and may retain stable spectral behavior under substantially weaker moment assumptions. Their high-dimensional spectral properties have consequently been studied in a variety of regimes; see, among others, \cite{bao2024spectral, BaiZhou,HeinyMikosch,HeinyComparison}. In the heavy-tailed setting, \cite{heinylimiting} identified a new family of limiting spectral distributions for sample correlation matrices, depending on the tail index and the aspect ratio. More recently, \cite{DongHeinyYao} obtained a resolvent-based description of the corresponding heavy-tailed Mar\v{c}enko--Pastur law through quadratic forms of self-normalized random vectors. Related spectral questions for correlation and self-normalized models have also been considered in \cite{HeinyParolya,DornemannFleermannHeiny}. These results make row-normalized random matrices natural objects in their own right, particularly when classical moment assumptions fail.

The present paper considers the non-Hermitian matrix itself rather than only its Gram matrix. This distinction is substantial. The spectrum of \(Y_nY_n^*\) records the singular values of \(Y_n\), but does not determine its complex eigenvalues. In particular, quantities such as the characteristic polynomial, traces, spectral radius, and the location of the complex spectrum contain genuinely non-Hermitian information that is invisible to the Gram matrix. This perspective is especially relevant in view of the recent work of Fran{\c c}ois, Heiny and Hu \cite{franccois2026characteristic}, who study the characteristic polynomial, trace fluctuations, and spectral radius of self-normalized non-Hermitian random matrices. Our results are complementary: we identify the limiting empirical eigenvalue distribution and study its geometry.

Let \(X_n=(X_{ij})_{1\le i,j\le n}\) be a matrix with i.i.d.\ real or complex entries, and define the row-normalized matrix \(Y_n=(Y_{ij})_{1\le i,j\le n}\) by
\[
    Y_{ij}
    =
    \frac{X_{ij}}
    {\left(\sum_{k=1}^n |X_{ik}|^2\right)^{1/2}} .
\]
We assume that the radial tail is regularly varying with index $\alpha\ge 0$, that is,
\[
    \mathbb P(|X_{11}|>t)
    =
    t^{-\alpha}L(t),
\]
where \(L\) is a slowly varying function (at infinity). The mechanism behind the model changes substantially with \(\alpha\).

Our main results concern the regime \(0<\alpha<2\). For every fixed \(z\in\mathbb C\), we prove, without any density assumption, convergence of the empirical singular value distribution of \(Y_n-zI\) to a deterministic probability measure. The limit depends only on \(\alpha\) and \(|z|\), and not on the slowly varying function or the angular tail distribution. Under the additional assumption that the entry distribution has a bounded density, we establish the logarithmic uniform integrability required by Hermitization and obtain convergence of the empirical eigenvalue distribution of \(Y_n\) to a deterministic probability measure \(\mu_\alpha^{\mathrm{SN}}\).

Our approach builds on the objective method and local weak convergence of Aldous and Steele \cite{aldous2004objective}, and more specifically on the Poisson weighted infinite tree (PWIT) and Hermitization framework developed for heavy-tailed non-Hermitian matrices in \cite{BCC}; see also \cite{bordenave2011reversible,bordenave2017spectrum}. After transposition, row normalization becomes column normalization. The raw heavy-tailed local limit can therefore be normalized column by column, leading to a directed, column-normalized PWIT. The dependence between weights belonging to the same column persists in the limiting operator and in the associated cavity equations. To pass from singular values to eigenvalues, we combine this local description with anti-concentration estimates for normalized rows and distance-to-subspace estimates adapted to arbitrary slowly varying tails. These estimates provide the small and intermediate singular value control required by Hermitization.

Closely related normalized PWIT structures already arise for heavy-tailed random Markov matrices. Bordenave, Caputo and Chafa\"i study the reversible setting in \cite{bordenave2011reversible}, where Poisson--Dirichlet weights and normalized PWIT-type structures appear, while Bordenave, Caputo, Chafa\"i and Piras treat the oriented setting in \cite{bordenave2017spectrum}. In the latter model, nonnegative heavy-tailed rows are normalized in \(L^1\), and the authors establish convergence of the shifted singular value distributions and, under a bounded density assumption, convergence of the empirical eigenvalue distribution to a rotationally invariant law supported on the unit disk. Thus the novelty of the present work lies not in the appearance of normalized PWITs themselves, but in the \(L^2\) normalization, which naturally covers the full range \(\alpha\in(0,2)\), in the treatment of arbitrary slowly varying tails and angular laws, and in the analysis of the resulting cavity equations, from which we derive detailed geometric information on the limiting eigenvalue distribution.

In particular, we prove that \(\mu_\alpha^{\mathrm{SN}}\) is rotationally invariant and supported in the closed unit disk, with outer support radius exactly one. Moreover, it has no point atoms and assigns no mass to the unit circle. A key ingredient is a general instability criterion for fresh  cavity recursions, which relates the outer support edge to the fractional moment behavior of the random weights appearing in the fresh side of the recursion. In the present \(L^2\)-normalized model these weights sum to one almost surely, leading to the exact outer radius one.

We formulate the support-edge argument as a general criterion for
measures admitting a fresh cavity recursion; see
Theorem~\ref{thm:general-support-criterion}.
Under its hypotheses, the squared outer support radius is
\[
    \inf_{0<p\le1}
       \left(\mathbb E\sum_k W_k^p\right)^{1/p}.
\]
where $W_k$ are the weights of the cavity recursion. The proof combines a fractional-moment upper bound with a supercriticality
argument for finite truncations of the recursion. 

We also identify the two parameter limits of our \(L^2\)-normalized model. As \(\alpha\uparrow2\), the measures \(\mu_\alpha^{\mathrm{SN}}\) converge to the circular law, whereas as \(\alpha\downarrow0\) they converge to the point mass at the origin. Thus self-normalization produces a family of non-Hermitian spectral laws interpolating between the circular law regime and an increasingly sparse one-hot regime, introduced in \cite{franccois2026characteristic}.

At zero shift, our singular value result gives an alternative operator description of the sample correlation problem studied in \cite{heinylimiting} and also yields a description of the moments of the limiting measure in terms of closed walks on the column-normalized PWIT. %The distinction between the two problems is important: the spectrum of \(Y_nY_n^*\) only requires the singular value law at \(z=0\), whereas the complex eigenvalue distribution of \(Y_n\) requires control of the full family of shifted matrices \(Y_n-zI\), \(z\in\mathbb C\). 
The normalized PWIT therefore provides a common framework for the Hermitian sample correlation problem and its non-Hermitian counterpart.

We finally treat the two endpoint mechanisms. In the normal domain of attraction, including the finite variance case, the row norms become asymptotically deterministic and the shifted singular value laws converge to their Gaussian counterparts; under a bounded density assumption, the circular law follows. At the opposite endpoint \(\alpha=0\), dominance of the maximum entry in a row implies that a typical normalized row is asymptotically concentrated on a single coordinate. The matrix is then close to the one-hot model introduced in \cite{franccois2026characteristic}; its empirical eigenvalue distribution converges to \(\delta_0\), while its squared singular value distribution converges to the Poisson law of mean one.

Finally, our results complement the spectral radius estimates of \cite{franccois2026characteristic}. For real entry distributions symmetric about zero, their results provide the upper bound one for the spectral radius. Under the bounded density assumption on \(X_{11}\), the fact that \(\mu_\alpha^{\mathrm{SN}}\) has outer support radius one gives the corresponding lower bound. Consequently, the spectral radius
\[
    \rho(Y_n)\longrightarrow1
\]
in probability for \(0<\alpha\le 2\). This universal behavior of the spectral radius contrasts with deterministically normalized heavy-tailed non-Hermitian matrices. Under the additional bounded density assumption, our results also answer the two open questions posed in \cite[Open Questions]{franccois2026characteristic}: the first is answered in the negative, while the second is answered in the affirmative.

\subsection{Notation}

Throughout the paper, $\mathbb{K} \in \{\mathbb{R},\mathbb{C}\}$, and we identify
$\mathbb{C}$ with $\mathbb{R}^2$. We write $d^2z$ for planar Lebesgue measure and 
\[
\overline{\mathbb{D}} := \{z \in \mathbb{C} : |z| \leq 1\}.
\]

For a matrix $A$, the notation $A^*$ denotes its conjugate transpose,
$\|A\|_{\mathrm{op}}$ its operator norm, and
$\|A\|_{\mathrm{F}}
:= \bigl(\operatorname{Tr} AA^*\bigr)^{1/2}$
its Frobenius norm. %If $A$ is an $n \times n$ matrix, then we write
%\[ s_1(A) \geq \cdots \geq s_n(A) \geq 0\]
%for its (ordered) singular values, and
%$\lambda_1(A),\ldots,\lambda_n(A)$ for its eigenvalues, counted with algebraic multiplicity.

For nonnegative sequences $(a_n)$ and $(b_n)$, we write
$a_n \sim b_n$ if
$a_n/b_n \longrightarrow 1$, as $\nto$.
The notation $a_n \lesssim b_n$ means that
$a_n \leq C b_n$
for a constant $C < \infty$ independent of $n$, while
$a_n \asymp b_n$ means that both $a_n \lesssim b_n$ and
$b_n \lesssim a_n$ hold. Dependence of the implicit constant on parameters is
indicated by a subscript; for example,
$a_n \lesssim_{\alpha} b_n$.

For a metric space $(E,d)$, let $\mathcal{P}(E)$ denote the set of Borel
probability measures on $E$. Weak convergence of probability measures is indicated by
$\mu_n \Rightarrow \mu$.

By $d_{\mathrm{BL}}(\mu,\nu)$ (resp.) $W_p(\mu,\nu)$ we denote the bounded Lipschitz distance (resp.) $p$-Wasserstein distance between two probability
measures $\mu$ and $\nu$. On a Polish space, $d_{\mathrm{BL}}$ metrizes weak convergence.

If $T:E\to F$ is measurable and $\mu\in\mathcal{P}(E)$, then
$T_{\#}\mu$ denotes the push-forward of $\mu$ under $T$.

We write
\[
\mathcal{D}(\mathbb{C})
:= C_c^\infty(\mathbb{C})
\simeq C_c^\infty(\mathbb{R}^2)
\]
for the space of smooth compactly supported test functions, and
$\mathcal{D}'(\mathbb{C})$ for its continuous dual, the space of distributions
on $\mathbb{C}$.
Finally, for a complex number $z=x+iy$, we use the Wirtinger derivatives
\[
\partial_z
:= \frac{1}{2}\bigl(\partial_x-i\partial_y\bigr),
\qquad
\partial_{\bar z}
:= \frac{1}{2}\bigl(\partial_x+i\partial_y\bigr),
\]
so that the Laplace differential operator $\Delta$ can be written as
$\Delta = 4\partial_z\partial_{\bar z}$.

\section{The model and main results}\label{sec:model-main-results}

Let $\mathbb K\in\{\R,\C\}$, with $\R$ viewed as a subset of $\C$,
and let $(X_{ij})_{i,j\ge1}$ be independent and identically distributed
$\mathbb K$-valued random variables with nondegenerate common law, that is, there is no $c\in \C$ such that $\mathbb{P}(X_{ij}=c)=1$.
Write $X$ for a random variable with this law and
$\overline F(t)=\Pp(|X|>t)$ for its radial tail.
Our heavy-tailed models are indexed by $\alpha\in[0,2]$ and satisfy
\begin{equation}\label{eq:model-radial-tail}
  \overline F(t)=t^{-\alpha}L(t),\qquad t>0,
\end{equation}
where $L$ is slowly varying, that is,   $\lim_{t\to \infty} \frac{L(ct)}{L(t)}=1$ for $c>0$.
For $\alpha=0$, this means that $\overline F$ itself is slowly varying;
in particular, $L(t)\to0$. At the endpoint $\alpha=2$, we also consider
the broader normal domain of attraction model, described in subsection \ref{section:normaldomain}, which includes
finite positive variance. 

For each $n$, put $X_n:=(X_{ij})_{1\le i,j\le n}$ and define
\begin{equation}\label{eq:model-row-normalization}
  D_{i,n}:=\left(\sum_{k=1}^n|X_{ik}|^2\right)^{1/2},
  \qquad
  Y_{ij}^{(n)}:=\frac{X_{ij}}{D_{i,n}},
  \qquad
  Y_n:=(Y_{ij}^{(n)})_{1\le i,j\le n}.
\end{equation}
If $D_{i,n}=0$, complete the $i$-th row $Y_n$ by putting a $1$ in a uniformly at random chosen position and setting the remaing entries $0$. If $p_0=\Pp(X=0)<1$, the probability that
this convention is used anywhere in $Y_n$ is at most $np_0^n$.
Since $\sum_n np_0^n<\infty$, this completion is almost surely irrelevant
for all sufficiently large $n$.
We usually suppress the superscript $(n)$. The rows of $Y_n$ are
independent, although their coordinates are generally dependent, and
\begin{equation}\label{eq:model-unit-rows}
  \sum_{j=1}^n|Y_{ij}|^2=1,
  \qquad \norm{Y_n}_F^2=n.
\end{equation}

For an $n\times n$ matrix $B$, let
$s_1(B)\ge\cdots\ge s_n(B)\ge0$ be its (ordered) singular values and let
$\lambda_1(B),\ldots,\lambda_n(B)$ be its eigenvalues, counted with
algebraic multiplicity. Define the measures 
\begin{equation}\label{eq:model-spectral-measures}
  \nu_B=\frac1n\sum_{j=1}^n\delta_{s_j(B)},\qquad
  H_B=\frac1n\sum_{j=1}^n\delta_{s_j(B)^2},\qquad
  \mu_B=\frac1n\sum_{j=1}^n\delta_{\lambda_j(B)}.
\end{equation}
Thus $H_B=\mu_{BB^*}$ is the empirical \emph{squared} singular value
measure, whereas $\mu_B$ is the complex eigenvalue measure.
The notation $B^*$ denotes conjugate transpose.

\subsection{The interior regime: \texorpdfstring{$0<\alpha<2$}{0 < alpha < 2}}
\label{ex:subsec-interior}

For the interior regime we will assume
\eqref{eq:model-radial-tail} with $0<\alpha<2$ and weak convergence
of the angular tail distribution: for some probability measure $\theta$
on the unit circle $\mathbb S^1:=\{w\in\mathbb C:|w|=1\}$,
\begin{equation}\label{eq:model-angular-tail}
  \mathcal L\left(\frac{X}{|X|}\,\middle|\,|X|>t\right)
  \Rightarrow\theta\,,\qquad t\to\infty.
\end{equation}
Equivalently, the corresponding probabilities $\P(X/|X|\in E ~|~ |X|>t)$ converge to $\theta(E)$
for every Borel set $E$ satisfying $\theta(\partial E)=0$.
For real entries, $\theta$ is supported on $\{-1,1\}$ and \eqref{eq:model-angular-tail} is simply a tail balance condition. For complex entries, $\theta$
may be any probability measure on the unit circle, including an atomic or
singular measure.
Condition \eqref{eq:model-angular-tail} is the weak angular form of
complex regular variation of (H1) in \cite{BCC}. 

For $n\ge2$, define the tail quantile
\begin{equation}\label{eq:model-tail-quantile}
  a_n=\max\{1,\inf\{a>0:\mathbb P(|X|>a)\le n^{-1}\}\},\qquad a_1=1.
\end{equation}
Then $a_n\to\infty$, $n\mathbb P(|X|>a_n)\to1$, and
$a_n=n^{1/\alpha}\ell(n)$ for a slowly varying $\ell$.
This auxiliary normalization is used only for the raw local limit, see Section \ref{sec:local}.
\medskip

For the convergence of the eigenvalue measure of $Y_n$, we will impose a bounded density assumption on
$\mathbb K$, 
\begin{equation}\label{eq:model-density}
  \mathbb P(X\in A)=\int_A f\,dm_{\mathbb K},
  \qquad \|f\|_{L^\infty(m_{\mathbb K})}<\infty.
\end{equation}
Here $m_{\mathbb R}$ is one-dimensional Lebesgue measure and
$m_{\mathbb C}$ is planar Lebesgue measure.

 %The limiting object is the column normalization of the directed, angularly marked Poisson weighted infinite tree (PWIT) arising in the paired-entry construction of Bordenave--Caputo--Chafa\"i \cite[Section~2.6]{BCC}, see Section \ref{sec:local} for more.
The following result establishes the limiting singular value and eigenvalue distributions of shifted matrices $Y_n$.
\begin{theorem}[singular value and eigenvalue limits]
\label{thm:interior-convergence}
Assume $0<\alpha<2$,
\eqref{eq:model-radial-tail} and
\eqref{eq:model-angular-tail}.
\begin{enumerate}[label=\textup{(\roman*)},leftmargin=*]
\item For every fixed $z\in\mathbb C$, there is a deterministic
probability measure $\nu_{\alpha,z}^{\mathrm{SN}}$ on $[0,\infty)$ such
that
\begin{equation}\label{eq:main-interior-limits}
  \nu_{Y_n-zI}\Rightarrow\nu_{\alpha,z}^{\mathrm{SN}}
  \qquad\text{in probability as }n\to\infty.
\end{equation}
 It depends only on $\alpha$ and $|z|$,
not on the slowly varying function or the angular law $\theta$.

\item Additionally if \eqref{eq:model-density} also holds true then for every fixed $z\in\mathbb C$, $\log$ is uniformly integrable
in probability for $(\nu_{Y_n-zI})_n$. More precisely, for every
$\varepsilon>0$,
\begin{equation}\label{ex:eq-main-ui}
 \lim_{\delta\downarrow0}\limsup_{n\to\infty}
 \mathbb P\left(
   \int_{[0,\delta]\cup[\delta^{-1},\infty)}
       |\log s|\,\nu_{Y_n-zI}(ds)>\varepsilon
 \right)=0,
\end{equation}
where $|\log0|=\infty$.

\item If additionally \eqref{eq:model-density} holds true,  then there is a deterministic probability measure
$\mu_\alpha^{\mathrm{SN}}$ on $\mathbb C$, depending only on $\alpha$,
such that
\begin{equation}\label{ex:eq-main-eigenvalues}
  \mu_{Y_n}\Rightarrow\mu_\alpha^{\mathrm{SN}}
  \qquad\text{in probability as }n\to\infty.
\end{equation}
For every fixed $z$, the integral
\[
  \Phi_\alpha(z)=\int_0^\infty\log s\,
                         \nu_{\alpha,z}^{\mathrm{SN}}(ds)
\]
is finite. This function is locally integrable, is the logarithmic
potential of $\mu_\alpha^{\mathrm{SN}}$ almost everywhere, and satisfies
\begin{equation}\label{eq:main-logarithmic-potential}
  \mu_\alpha^{\mathrm{SN}}
      =\frac1{2\pi}\Delta\Phi_\alpha
      \qquad\text{in }\mathcal D'(\mathbb C).
\end{equation}
\end{enumerate}
\end{theorem}

At $z=0$, pushing the singular value limit forward by $s\mapsto s^2$
gives the limiting eigenvalue measure of $Y_nY_n^*$. Hence, $\nu_{\alpha,0}^{\mathrm{SN}}$ is the $\alpha$-heavy Marchenko-Pastur law with parameter $1$, which was identified through its moments in \cite{heinylimiting}. As a by-product, our method yields an alternative description of those moments in terms of closed walks on a column normalized PWIT.

\subsection{Properties of the limiting eigenvalue measure}

\begin{figure}\label{figure}
    \centering
    \includegraphics[width=0.95\linewidth]{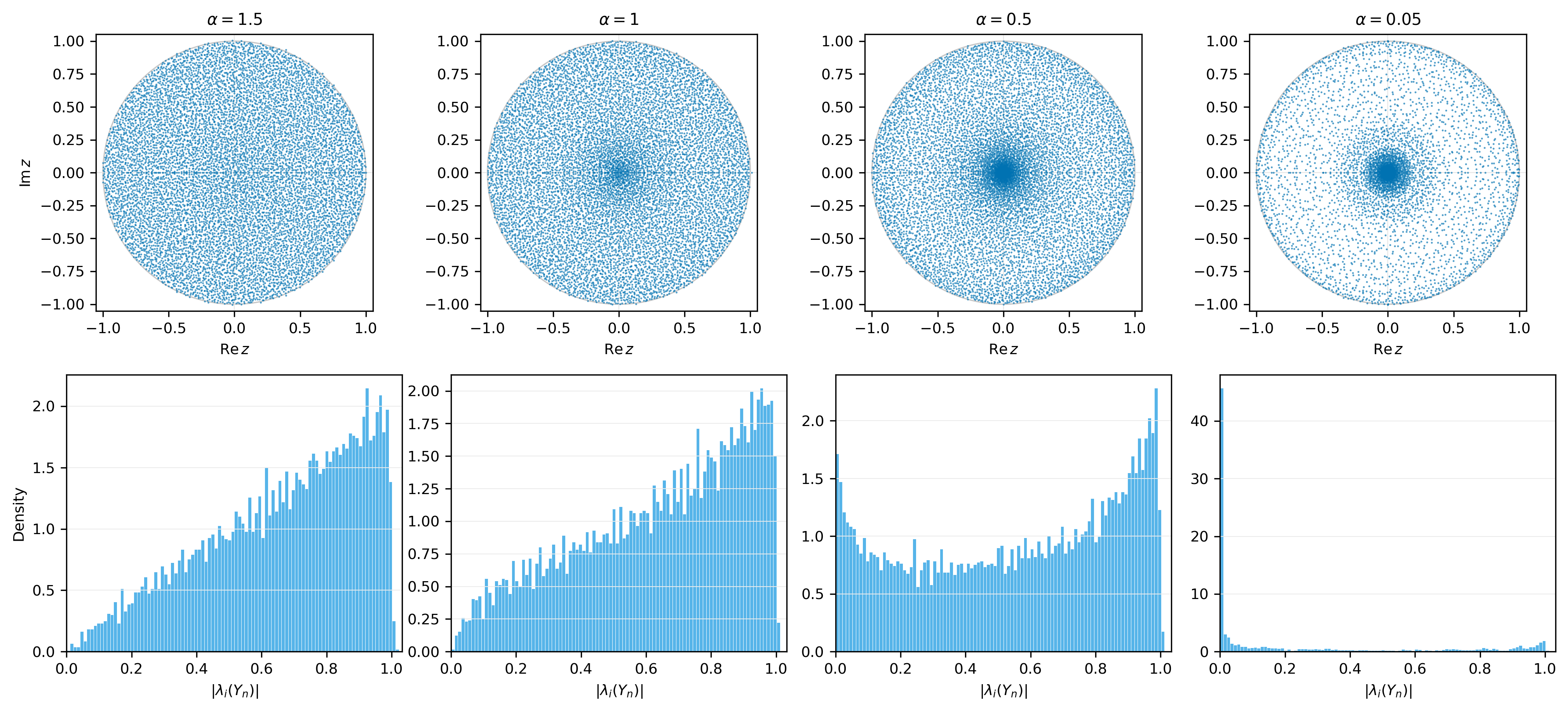}
    \caption{Top row: Eigenvalues of $Y_n$ for different values of $\alpha\in\{1.5,1,0.5,0.05\}$ with $n=10000$.  The entries are symmetric Pareto variables satisfying $\P(|X_{ij}|>x)=x^{-\alpha}$ for $x\ge1$. Bottom row: Corresponding histograms of the moduli of the eigenvalues.}
    \label{fig:eigenvalue-clouds}
\end{figure}

Next, we present poperties of $\mu_\alpha^{\mathrm{SN}}$ which will be illustrated in Figure~\ref{fig:eigenvalue-clouds}.
\begin{theorem}[Geometry of the limiting measure]
\label{thm:limit-properties}
For $\alpha\in (0,2)$, the
universal measure $\mu_\alpha^{\mathrm{SN}}$ has the following properties.
\begin{enumerate}[label=\textup{(\roman*)},leftmargin=*]
\item It is supported in the closed unit disk and is rotationally
invariant about the origin (isotropic):
\[
  \operatorname{supp}\mu_\alpha^{\mathrm{SN}}
  \subseteq\overline{\mathbb D},\qquad
  \mu_\alpha^{\mathrm{SN}}(e^{i\vartheta}E)
  =\mu_\alpha^{\mathrm{SN}}(E)
\]
for every Borel set $E\subseteq\C$ and every $\vartheta\in\R$.

\item Its outer support radius is exactly one:
\[
  \sup\{|w|:w\in\operatorname{supp}\mu_\alpha^{\mathrm{SN}}\}=1.
\]
Equivalently,
$\mu_\alpha^{\mathrm{SN}}\{|w|>r\}>0$ for every $0\le r<1$. 

\item It has no point atoms:
\[
  \mu_\alpha^{\mathrm{SN}}(\{w\})=0,
  \qquad w\in\C.
\]

\item It assigns no mass to the unit circle:
\[
  \mu_\alpha^{\mathrm{SN}}\{w:|w|=1\}=0.
\]
\end{enumerate}
\end{theorem}

We also include the following convergence theorem for the parameter endpoints.
\begin{theorem}[Parameter endpoints]\label{thm:parameterendpoints}

For the collection of measures $(\mu_{\alpha}^{\mathrm{SN}})_{\alpha\in (0,2)}$ we have the following limits.
   \begin{enumerate}[label=\textup{(\roman*)},leftmargin=*]
       \item At the upper endpoint, the circular law appears:
\[
  \mu_\alpha^{\mathrm{SN}}\Rightarrow\mu_{\mathrm{circ}}
  \qquad\text{as }\alpha\uparrow2.
\]

\item At the lower endpoint, the Dirac mass in zero appears:
\[
  \mu_\alpha^{\mathrm{SN}}\Rightarrow\delta_0
  \qquad\text{as }\alpha\downarrow0.
\]
   \end{enumerate} 

\end{theorem}
 The two parameter limits above concern the family of
limiting measures: $n\to\infty$ is taken first at fixed
$\alpha\in(0,2)$. The exact endpoint models are addressed separately
below.  

Figure~\ref{fig:eigenvalue-clouds} illustrates the conclusions of
Theorems~\ref{thm:limit-properties}(ii)
and~\ref{thm:parameterendpoints}(ii) for $n=10\,000$ and symmetric
Pareto entries with $\alpha\in\{1.5,1,0.5,0.05\}$.
The spectrum extends towards the boundary of the unit disk,
while the eigenvalues become increasingly concentrated near
the origin as $\alpha$ decreases.

\subsection{
  \texorpdfstring{Normal domain of attraction}{alpha = 2}}\label{section:normaldomain}

At the endpoint $\alpha=2$, and in general in the normal domain of attraction, the stable-PWIT description used for
$0<\alpha<2$ is replaced by relative stability of the Euclidean row
norms. 
Let $\mathbb K\in\{\mathbb R,\mathbb C\}$, and suppose that the entries
are i.i.d.\ copies of a nondegenerate $\mathbb K$-valued random
variable $X$ satisfying
\[
  \mathbb E X=0.
\]
An equivalent condition for normal domain of attraction is that the truncated second
moment
\begin{equation}\label{eq:model-truncated-second-moment}
  \ell_2(x)
  :=
  \mathbb E\bigl[
    |X|^2\mathbf1_{\{|X|\le x\}}
  \bigr]
\end{equation}
is slowly varying at infinity.
This condition includes every centered distribution with finite
positive second absolute moment, in which case $\ell_2(x)\longrightarrow\mathbb E|X|^2$.
It also includes regularly varying radial tails of index two, that is,
$\mathbb P(|X|>x)=x^{-2}L(x)$
with $L$ slowly varying. In the latter case, $\ell_2$ may diverge
slowly. 
\medskip

Now, let $G_n$ be and $n\times n$ matrix with i.i.d.\  standard Gaussian entries over $\mathbb K$, normalized
so that $\mathbb E|G_{11}|^2=1$.
Thus, in the real case, $G_{11}\sim N(0,1)$, while in the complex case
$G_{11}$ is a standard circular complex Gaussian. Denote by
$\nu_z^{\mathrm G_n}$ the singular value measure
of $\frac1{\sqrt n}G_n-zI$.
Then there exists a deterministic probability measure $\nu_z^{\mathrm{G}}$ such that $\nu_z^{\mathrm{G}_n}\Rightarrow\nu_z^{\mathrm G}$ in probability, see Theorem~4.9 of
\cite{bordenave2012around}. 
The next theorem shows that in the normal domain of attraction,  the (shifted) singular value and eigenvalues measures of $Y_n$ are universal in the sense that they have the same limits as in the Gaussian case.

\begin{theorem}[The normal-domain circular law]
\label{thm:exact-two}
Let $\mathbb K\in\{\mathbb R,\mathbb C\}$, and let the entries be
i.i.d.\ copies of a nondegenerate $\mathbb K$-valued random variable
$X$ satisfying $\mathbb E X=0$.
Assume that the truncated second moment
\eqref{eq:model-truncated-second-moment} is
slowly varying. Then, for every fixed $z\in\mathbb C$,
\begin{equation}\label{eq:exact-two-singular-values}
  \nu_{Y_n-zI}\Rightarrow\nu_z^{\mathrm G}
  \qquad\text{in probability as }n\to\infty.
\end{equation}

If, in addition, the law of $X$ has a bounded density with respect to
Lebesgue measure on $\mathbb K$, then
\begin{equation}\label{eq:exact-two-circular-law}
  \mu_{Y_n}\Rightarrow\mu_{\mathrm{circ}}
  \qquad\text{in probability as }n\to\infty.
\end{equation}
\end{theorem}
It is worth mentioning that, at zero shift, the limiting singular value measure is the
quarter-circle law, that is,
\[
  \nu_0^{\mathrm G}(ds)
  =
  \frac1\pi\sqrt{4-s^2}\,
  \mathbf1_{\{0<s<2\}}\,ds.
\]
The proof of Theorem~\ref{thm:exact-two} is inspired by \cite{bordenave2012circular}, where convergence to the circular law for random Markov matrices whose weights have finite variance was established.

\subsection{Spectral radius for $\alpha \in (0,2]$}
Furthermore, we present limit theory for the spectral radius of $Y_n$ for all positive $\alpha>0$ under the extra hypothesis for $X_{1,1}$ to be real and symmetric. In all cases for $\alpha$ the spectral radius converges to $1$. We prove this, using the upper bound provided by \cite{franccois2026characteristic} and the lower bound provided by Theorem \ref{thm:exact-two} and Theorem \ref{thm:limit-properties}(ii).

\begin{theorem}[Spectral radius for real symmetric entries]
\label{cor:real-symmetric-radius}
Assume the hypotheses of Theorem~\ref{thm:interior-convergence} including \eqref{eq:model-density} or \ref{thm:exact-two}, and
suppose additionally that $\mathbb K=\R$ and $X\stackrel d=-X$.
Then
\[
  \rho(Y_n):=\max_{1\le j\le n}|\lambda_j(Y_n)|
  \longrightarrow1
  \qquad\text{in probability}.
\]
\end{theorem}
\begin{proof}
We will give an upper and lower bound for $\rho(Y_n)$.

\begin{itemize}
    \item Upper bound:
Under the assumptions of Theorem \ref{thm:interior-convergence}, for every $\varepsilon>0$,
\cite[Corollary~2.6]{franccois2026characteristic} gives
\[
 \mathbb P\{\rho(Y_n)\ge1+\varepsilon\}\longrightarrow0.
\]

Similarly under the assumptions of Theorem \ref{thm:exact-two} and using \cite[Theorem 2.12]{franccois2026characteristic} one can conclude accordingly.
\item Lower bound: From Theorem \ref{thm:limit-properties} (ii) respectively Theorem \ref{thm:exact-two}, one immediately  deduces that for all $\epsilon>0$ it holds
\[
 \mathbb P\{\rho(Y_n)\le 1-\varepsilon\}\longrightarrow0.
\]  
\end{itemize}
\end{proof}

A similar spectral radius limit holds for random matrices with i.i.d.\ entries of mean zero and variance one: after the deterministic normalization $n^{-1/2}$, the spectral radius converges in probability to $1$, the outer radius of the support of the circular law; see \cite[Theorem 1.1]{bordenave2022convergence}. In contrast, for i.i.d. entries satisfying \eqref{eq:model-radial-tail} with $0<\alpha<2$, and under the additional assumptions of \cite{BCC}, the empirical eigenvalue distribution of the deterministically normalized matrix $a_n^{-1}X_n$ converges almost surely to a probability measure with unbounded support; see \cite[Theorems 1.2 and 1.3]{BCC}. Consequently, $\rho(a_n^{-1}X_n)\to\infty$ almost surely. See also \cite[Section 1.1.7]{bordenave2022convergence} for a discussion of this divergence and the associated question of the limiting point process at the spectral edge.
This comparison indicates that $\rho(Y_n)$ exhibits a remarkable robustness. 

\subsection{The slowly varying endpoint: \texorpdfstring{$\alpha=0$}{alpha = 0}}

At $\alpha=0$, one coordinate carries asymptotically all the squared mass of a typical row of $Y_n$. The next theorem shows how this sparsity affects the limiting eigenvalue distributions.
\begin{theorem}[Degenerate limit at $\alpha=0$]
\label{thm:exact-zero}
Let the entries of $X_n$ be real-valued, i.i.d., nondegenerate random variables
whose radial tail $\Pp(|X|>t)$ is slowly varying.
Then
\[
  \mu_{Y_n}\Rightarrow\delta_0
  \qquad\text{almost surely as }n\to\infty.
\]
\end{theorem}
Theorem~\ref{thm:exact-zero} reveals a stark contrast between the eigenvalue and singular value distributions of $Y_n$. Recently, \cite{DongHeinyYao} proved that $\nu_{Y_nY_n^*}$ converges to
%\[  \Pi_1=e^{-1}\sum_{k=0}^\infty\frac1{k!}\delta_k,\]
the Poisson distribution of mean one. In contrast, our Theorem~\ref{thm:exact-zero} shows that $\mu_{Y_n}$ becomes degenerate.
%Equivalently, the singular values themselves satisfy
%\[
 % \nu_{Y_n}\Rightarrow\nu_0^{\mathrm{SN}}
  %:=e^{-1}\sum_{k=0}^\infty\frac1{k!}\delta_{\sqrt{k}}
  %\qquad\text{almost surely}.
%\]

\section{Proof of Theorem~\ref{thm:interior-convergence}}
\label{sec:local}

{\it Throughout this section $0<\alpha<2$ is fixed, the radial tail is
regularly varying with index $\alpha$, and angular convergence means
\eqref{eq:model-angular-tail}. The density hypothesis is only used for
the least singular value estimate, Proposition \ref{ex:prop-invertibility}. Constants may depend on the entry
law and on $\alpha$. The slowly varying function in the radial tail is assumed to be fixed. We use ordinary matrix
coefficients, so $a_{uv}$ means the $u$-th coordinate of $\mathcal A e_v$.
The notation $\sharp$ denotes the formal adjoint on finitely supported
inputs.}

The proof follows three steps. We construct the angularly marked PWIT,
normalize its columns, and obtain the shifted singular value limit by
root-resolvent convergence. We then prove polynomial invertibility and
fractional inverse-distance estimates for general regularly varying
real or complex entries, giving logarithmic uniform integrability.
Finally, a local integrability argument passes from logarithmic
potentials to the eigenvalue measure. Phase removal and rotation of the
shift yield the asserted universality.

\subsection{The raw PWIT and column normalization}
\label{ex:subsec-pwit}

Put
\[
  \widehat A_n:=a_n^{-1}X_n^{\mathsf T} \qquad \text{ and } \qquad M_n:=Y_n^{\mathsf T}.
\]
For every $z\in\mathbb C$ and every $j$, we have
\begin{equation}\label{ex:eq-transpose}
 \det(M_n-zI)=\det(Y_n-zI),\qquad
 s_j(M_n-zI)=s_j(Y_n-zI),
\end{equation}
where for the second identity we used that the squared singular values of a matrix  $C^{\mathsf T}$
are the eigenvalues of $\overline{CC^*}$ and hence those of $CC^*$. We work with the transpose $M_n$ because the ``local weak convergence'' of the raw model $\widehat A_n$, \cite[Theorem 2.10]{BCC}, can be extended to the local weak convergence of the column normalized matrices $M_n$; see Proposition \ref{ex:prop-local-normalized}.

If $d_{n,j}=\|\widehat A_ne_j\|_2$, then
\begin{equation}\label{ex:eq-finite-normalization}
 M_ne_j=\frac{\widehat A_ne_j}{d_{n,j}}
\end{equation}
and the denominators are eventually almost surely positive as discussed below \eqref{eq:model-row-normalization}. Moreover, the structure of $M_n$ is preserved under multiplication with permutation matrices.

We label the vertices of the rooted tree by finite sequences of positive integers. With the convention $\mathbb N^0:=\{\varnothing\}$, we set
\[
V:=\bigcup_{h\ge0}\mathbb N^h,
\qquad o:=\varnothing.
\]
Thus, $o$ is the root, and a vertex $v=(v_1,\ldots,v_h)\in\mathbb N^h$ lies at generation $h$. For $k\in\mathbb N$, we write $vk:=(v_1,\ldots,v_h,k)$ for its $k$th child. We call $v$ the parent of $vk$.

Let $\ell^2(V)$ be the Hilbert space of square-summable complex-valued functions on $V$, with its canonical orthonormal basis $(e_v)_{v\in V}$, where
\[
e_v(w):=\mathbf 1_{\{w=v\}},
\qquad v,w\in V.
\]
We denote by
\[
\mathcal D(V):=\operatorname{span}\{e_v:v\in V\}
\]
the dense subspace of finitely supported functions on $V$.
Independently at every vertex $v$, take the arrival times
$\Gamma_{v,1}<\Gamma_{v,2}<\cdots$ of a rate-$2$ Poisson process.
Furthermore, for every $v,k$ consider
independent phases $\Theta_{v,k}$ with common law $\theta$ and
independent marks $\varepsilon_{v,k}\sim\operatorname{Bernoulli}(1/2)$.
The families $\{\Gamma_{v,k}\}$, $\{\Theta_{v,k}\}$ and $\{\varepsilon_{v,k}\}$ are mutually independent. Next, we define the operator $\mathcal A:\mathcal D(V)\to\ell^2(V)$ by
\[
\mathcal A e_v:=\sum_{u\in V}a_{u,v}e_u,
\qquad v\in V,
\]
where, for each $v\in V$ and $k\in\mathbb N$, we set
\[
\zeta_{v,k}:=\Theta_{v,k}\Gamma_{v,k}^{-1/\alpha},
\qquad
a_{v,vk}:=\varepsilon_{v,k}\zeta_{v,k},
\qquad
a_{vk,v}:=(1-\varepsilon_{v,k})\zeta_{v,k}.
\]
All other coefficients, including the diagonal coefficients, are set equal to zero. Define also
\[
 \mathcal A^\sharp e_v:=\sum_u\overline{a_{vu}}e_u.
\]
Since $\Gamma_{v,k}/k\to1/2$ almost surely as $k\to \infty$ and $2/\alpha>1$, we get $\sum_{k\ge1}|\zeta_{v,k}|^2<\infty$
simultaneously for all $v$. There is at most one additional parent
coefficient in a row or column. Thus, the operators operators $\mathcal{A}$ and $\mathcal A^\sharp$ are
well-defined on $\mathcal D(V)$.

In the following, we write $\beta:=\alpha/2$ so that the squared child magnitudes
$J_{v,k}:=|\zeta_{v,k}|^2$ have Poisson intensity
$2\beta x^{-1-\beta}\,dx$. Each orientation-thinned cloud has intensity
\begin{equation}\label{ex:eq-squared-intensity}
 \Lambda_\beta(dx):=\beta x^{-1-\beta}\,dx.
\end{equation}
If $J_v,\varepsilon_v$ are the parent-edge marks at $v\ne o$, then
\begin{equation}\label{ex:eq-normalizer}
 d(v)^2:=\sum_u|a_{uv}|^2
 =\mathbf1_{\{v\ne o\}}\varepsilon_vJ_v
       +\sum_k(1-\varepsilon_{v,k})J_{v,k}.
\end{equation}
The fresh sum is finite and strictly positive almost surely (see Lemma \ref{ex:lem-reverse}). In
particular $0<d(v)<\infty$ at every vertex. Define
\begin{equation}\label{ex:eq-M-definition}
 \mathcal M e_v:=\sum_u m_{uv}e_u \qquad \text{ with } \quad m_{uv}:=a_{uv}/d(v)
\end{equation}
and note that $\sum_u|m_{uv}|^2=1$.
The parent contribution in \eqref{ex:eq-normalizer} is retained when
an edge is cut from a cavity operator. The normalizers at distinct
vertices are not asserted to be independent.
\subsection{Properties of the limiting operator}
\label{ex:subsec-esa}

\begin{lemma}[Reverse square summability]\label{ex:lem-reverse}
Almost surely $\sum_u|m_{vu}|^2<\infty$ for every $v$. At the root,
\begin{equation}\label{ex:eq-root-row-mean}
 \mathbb E\sum_u|m_{ou}|^2=1.
\end{equation}
\end{lemma}

\begin{proof}
For a child contributing to row $v$, the connecting squared magnitude
$j$ belongs to the child's column. Its squared normalizer is $j+S$,
where $S$ is an independent fresh source sum with intensity
\eqref{ex:eq-squared-intensity}. These fresh marks are independent
across children. Consequently,
\[
 \sum_u|m_{vu}|^2
 \le \mathbf1_{\{v\ne o\}}
    +\sum_{k:\varepsilon_{v,k}=1}\frac{J_{v,k}}{J_{v,k}+S_k}.
\]
The Poisson Laplace functional gives
\[
 \mathbb Ee^{-uS}=e^{-\Gamma(1-\beta) \, u^\beta},\qquad
 \mathbb ES^{-\beta}
    =\frac1{\Gamma(1+\beta)\Gamma(1-\beta)},
\]
where the second identity follows by writing
$s^{-\beta}=\Gamma(\beta)^{-1}\int_0^\infty
 u^{\beta-1}e^{-us}\,du$ and using Tonelli.
Applying Campbell's formula to the independently marked point process, see \cite[Propositions 2.7 and 5.5]{LastPenrose}, and then using Tonelli's theorem and the substitution $x=Sy$, we obtain
\begin{align*}
 \mathbb E\sum_{k:\varepsilon_{v,k}=1}
      \frac{J_{v,k}}{J_{v,k}+S_k}
 &=\beta\int_0^\infty x^{-1-\beta}
                    \mathbb E\frac{x}{x+S}\,dx\\
 &=\beta\mathbb ES^{-\beta}
        \int_0^\infty\frac{y^{-\beta}}{1+y}\,dy\\
  &=\beta\mathbb ES^{-\beta}\,
        \Gamma(\beta)\Gamma(1-\beta)      =1.
\end{align*}
%The last integral is $\Gamma(\beta)\Gamma(1-\beta)$.
Finite expectation of this nonnegative series gives almost-sure
finiteness. The root has no parent term, so its expectation is exactly
one. Countability proves the claim for all $v \in V$.
\end{proof}

Let $V^b:=V^+\sqcup V^-$, with basis $e_v^+,e_v^-$ as in \cite{BCC}. On
$\mathcal D(V^b)$ define the minimal symmetric operator $\mathcal B_0(0)$
by coefficient $m_{uv}$ from $v^-$ to $u^+$ and coefficient
$\overline{m_{uv}}$ in the reverse direction, and zero otherwise.
Both rows and columns are square summable, so the operator is
well-defined. Write
\[
 \mathcal B_0(0)=
 \begin{pmatrix}0&\mathcal M\\\mathcal M^\sharp&0\end{pmatrix}
 \quad\text{on }\mathcal D(V^b).
\]

\begin{lemma}[One-sided $\ell^2$ exhaustion criterion]
\label{ex:lem-exhaustion}
Let $G=(W,F)$ be a countable bipartite tree, with
$W=W_0\sqcup W_1$, and let its edge weights satisfy
$b_{xy}=\overline{b_{yx}}$ whenever $\{x,y\}\in F$.
Assume that
\[
  \sup_{y\in W_1}\sum_{x\sim y}|b_{xy}|^2\le C<\infty
  \qquad \text{ and } \qquad 
  \sum_{y\sim x}|b_{xy}|^2<\infty
  \quad\text{for every }x\in W_0.
\]
and define the weighted adjacency operator $T_0$ on $\ell^2(W)$ by 
\[
  D(T_0)=\mathcal D_{\mathrm{fin}}(W)
  :=\{f\in\ell^2(W):\operatorname{supp}f\text{ is finite}\},
  \qquad
  (T_0f)(x):=\sum_{y\sim x}b_{xy}f(y).
\]
Then $T_0$ is essentially self-adjoint. In particular, its closure is
self-adjoint and $\mathcal D_{\mathrm{fin}}(W)$ is a core for that
closure.
\end{lemma}

\begin{proof}
The hypotheses and the Hermitian symmetry imply square summability of
all rows and columns of $(b_{xy})$. Consequently $T_0e_x\in\ell^2(W)$
for every coordinate vector $e_x$, so $T_0$ is well-defined on
$\mathcal D_{\mathrm{fin}}(W)$. This domain is dense, and
$b_{xy}=\overline{b_{yx}}$ implies that $T_0$ is symmetric.
If $W$ is finite, the assertion is immediate. We therefore assume that
$W$ is countably infinite.

\medskip
Set $\kappa=\max\{C,1\}$. For every $x\in W_0$, choose a finite set
$F_x\subseteq\{y:y\sim x\}$ such that
\[
  \sum_{\substack{y\sim x\\y\notin F_x}}|b_{xy}|^2\le\kappa.
\]
This is possible by square summability at $x$.
Enumerate $W=\{w_1,w_2,\ldots\}$ and set $S_0=\varnothing$.
Given $S_{k-1}$, let $H_k$ be the connected hull of
$S_{k-1}\cup\{w_k\}$ in $G$, and put
\[
  S_k=H_k\cup\bigcup_{x\in H_k\cap W_0}F_x.
\]
The set $H_k$ is finite: it is obtained by joining finitely many
vertices by the unique finite paths between them. Thus $S_k$ is also
finite and connected. Moreover,
  $S_{k-1}\subseteq S_k$ and $\bigcup_{k\ge1}S_k=W$.

We claim that
\begin{equation}\label{ex:eq-exhaustion-boundary}
  \sum_{\substack{y\notin S_k\\y\sim x}}|b_{xy}|^2\le\kappa
  \qquad\text{for every }k\ge1\text{ and }x\in S_k.
\end{equation}
If $x\in S_k\cap W_1$, its entire incident squared weight is at most
$C\le\kappa$. If $x\in S_k\cap W_0$, then $x\in H_k$, because all
vertices added to $H_k$ belong to $W_1$. Hence $F_x\subseteq S_k$,
and the choice of $F_x$ gives \eqref{ex:eq-exhaustion-boundary}. The claim now follows by Corollary 2.4. of \cite{BCC}.
\end{proof}

\begin{proposition}[Self-adjoint realization and common core]
\label{ex:prop-esa}
Almost surely $\mathcal B_0(0)$ is essentially self-adjoint. Denote its
closure by $\mathcal B(0)$. For every $z\in\mathbb C$, let
\[
 S(z)e_v^-=ze_v^+,\qquad S(z)e_v^+=\overline z e_v^-,
 \qquad \mathcal B(z)=\mathcal B(0)-S(z).
\]
Then $\mathcal B(z)$ is self-adjoint on $D(\mathcal B(0))$, and
$\mathcal D(V^b)$ is a core for all these operators.
\end{proposition}

\begin{proof}
The claim follows by Lemmas \ref{ex:lem-exhaustion} and \ref{ex:lem-reverse}.
\end{proof}
\subsection{Weak convergence to the operator}
We identify $[n]$ with the one-letter words
$\{(1),\ldots,(n)\}\subset V$ and regard $\widehat A_n$ as an
operator on $\ell^2(V)$ by extending its matrix coefficients
by zero outside $[n]\times[n]$.

\begin{lemma}[Definition 2.6 and the proof of Theorem 2.10 in \cite{BCC}]
\label{ex:lem-raw-local}
On a suitable coupling, there exist random bijections
$\sigma_n:V\to V$, with $\sigma_n(o)=1$, and associated
permutation unitary operators $P_n$ on $\ell^2(V)$, defined by
\[
    P_n e_v=e_{\sigma_n(v)},\qquad v\in V,
\]
such that, setting
\[
    V_n:=\sigma_n^{-1}([n]),
    \qquad
    T_n:=P_n^{-1}\widehat A_nP_n,
\]
one has, almost surely, simultaneously for every $v\in V$,
\begin{equation}\label{ex:eq-raw-joint}
    T_ne_v\longrightarrow\mathcal A e_v,
    \qquad
    T_n^*e_v\longrightarrow\mathcal A^\sharp e_v
    \quad\text{in }\ell^2(V).
\end{equation}
Moreover, $|V_n|=n$, $o\in V_n$, and almost surely every
finite subset of $V$ is contained in $V_n$ for all
sufficiently large $n$.
\end{lemma}
\begin{proposition}[Column convergence after normalization]
\label{ex:prop-local-normalized}
Under the coupling of Lemma~\ref{ex:lem-raw-local},
\[
 P_n^{-1}M_nP_ne_v\longrightarrow\mathcal M e_v
       \quad\text{in }\ell^2(V)
\]
almost surely for each fixed $v$, and hence on every finitely supported
input vector.
\end{proposition}

\begin{proof}
For each fixed $v$, eventually its relabeled column belongs to the
finite matrix and is nonzero. Permutation equivariance and
\eqref{ex:eq-finite-normalization} give
\[
 P_n^{-1}M_nP_ne_v=\frac{T_ne_v}{\|T_ne_v\|_2}.
\]
The map $x\mapsto x/\|x\|_2$ is continuous at
$\mathcal A e_v\ne0$. Apply Lemma~\ref{ex:lem-raw-local}, then take
finite linear combinations.
\end{proof}

\subsection{Root resolvents and the shifted singular value limit}
\label{ex:subsec-resolvents}

For each $n \in \N$ and $z\in \C$ define the matrix $B_n(z)\in \mathcal{M}_{n}(\mathcal{M}_2(\C))\simeq \mathcal{M}_{2n}(\C)$ with entries
\[
 (B_n(z))_{i,j}=\begin{pmatrix}0&(M_n-zI)_{i,j}\\(M_n^*-\overline zI)_{i,j}&0\end{pmatrix},
 \qquad i,j=1,\ldots,n\,.
\]
Moreover, let
$$G_{n,j}(z,\eta)=\bigl[(B_n(z)-\eta I)^{-1}\bigr]_{jj},$$
where $[\,\cdot\,]_{ij}$ indicates the $2\times2$ block on the two copies
of vertices $i$ and $j$.

For each \(v\in V\), let \(\Pi_v:\ell^2(V^{\mathrm b})\to\mathbb C^2\) denote the coordinate projection onto the two copies \(v^+\) and \(v^-\) of \(v\), defined by
\[
\Pi_v f=
\begin{pmatrix}
f(v^+)\\
f(v^-)
\end{pmatrix}.
\]
Its adjoint is the canonical embedding
\[
\Pi_v^*
\begin{pmatrix}
a\\ b
\end{pmatrix}
=
a e_v^+ + b e_v^-.
\]
For \(z\in\mathbb C\) and \(\eta\in\mathbb C_+\), we define the resolvent block at \(v\) by
\[
G_v(z,\eta)
:=
\Pi_v(\mathcal B(z)-\eta I)^{-1}\Pi_v^*
\in\mathbb C^{2\times2},
\qquad v\in V.
\]
Thus \(G_v(z,\eta)\) is the compression of the resolvent to
\(\operatorname{span}\{e_v^+,e_v^-\}\), expressed in this ordered basis. In particular, \(G_o(z,\eta)\) denotes the resolvent block at the root.
Throughout this section, $\eta\in\mathbb C_+$. 
\begin{proposition}[Root-block convergence]\label{ex:prop-root}
For fixed $z\in\mathbb C$ and $\eta\in\mathbb C_+$,
\[
 G_{n,1}(z,\eta)\Rightarrow G_o(z,\eta)\,, \qquad \nto.
\]
\end{proposition}
\begin{proof}
    The claim follows by Proposition \ref{ex:prop-local-normalized}, Proposition \ref{ex:prop-esa} and Theorem 2.7 in \cite{BCC}.
\end{proof}

For a self-adjoint operator $T$ and unit vector $e$, write
$\rho_{T,e}(E)=\langle e,\mathbf1_E(T)e\rangle$. Define
\begin{equation}\label{ex:eq-limit-sv-definition}
 \check\nu_{\alpha,\theta;z}
   =\frac12\mathbb E\bigl[
       \rho_{\mathcal B(z),e_o^+}+\rho_{\mathcal B(z),e_o^-}\bigr],
 \qquad
 \nu_{\alpha,\theta;z}=|\cdot|_\#\check\nu_{\alpha,\theta;z}.
\end{equation}
Recall that the column-normalized Poisson weighted infinite tree, introduced at the beginning of Section \ref{sec:local}, depends on the angular measure $\theta$ appearing in \eqref{eq:model-angular-tail}.

Both are probability measures on $\R$. The first is symmetric because the
sign-change unitary on the minus side conjugates $\mathcal B(z)$
to $-\mathcal B(z)$ while fixing either root vector up to sign.
Its Stieltjes transform is
\begin{equation}\label{ex:eq-root-stieltjes}
 m_{\check\nu_{\alpha,\theta;z}}(\eta)
   =\frac12\mathbb E\operatorname{Tr}G_o(z,\eta).
\end{equation}

\begin{proposition}[Weak shifted singular value convergence]
\label{ex:prop-weak-sv}
For each fixed $z \in \C$, it holds $\nu_{Y_n-zI}\Rightarrow\nu_{\alpha,\theta;z}$
in probability, as $\nto$.
\end{proposition}

\begin{proof}
The empirical measure of $B_n(z)$ is the symmetrization
$\check\nu_{n,z}$ of $\nu_{M_n-zI}$, since its eigenvalues are the
numbers $\pm s_j(M_n-zI)$; see \cite[Theorem~2.1]{BCC}. Put
\[
 m_{n,z}(\eta)=\frac1{2n}\operatorname{Tr}(B_n(z)-\eta I)^{-1}.
\]
Permutation invariance of the law of the diagonal coordinates of the resolvent, Proposition~\ref{ex:prop-root}, and the bounded
resolvent entries give
\begin{align}\label{eq-convergence-of-mean-shifted-resolvent}
 \mathbb E m_{n,z}(\eta)
  =\tfrac12\mathbb E\operatorname{Tr}G_{n,1}(z,\eta)
  \longrightarrow m_{\check\nu_{\alpha,\theta;z}}(\eta).
\end{align}
Since the columns of $M_n$ are independent, an application of \cite[Lemma C.4]{bordenave2013localization} yields that 
\begin{align}\label{eq-concentration-shifted-stieltjes}
    \mathbb{P}\left( \left| m_{n,z}(\eta)-\mathbb{E} m_{n,z}(\eta) \right|\geq t\right)\leq 2 \exp(-c n \Im^2(\eta) t^2)
\end{align}
for some absolute constant $c>0$. Thus, the Stieltjes transform concentrates around its mean asymptotically a.s. and so the claim follows by \eqref{eq-convergence-of-mean-shifted-resolvent} and \eqref{eq-concentration-shifted-stieltjes}
\end{proof}

\subsection{Removing the phases and rotating the shift}
\label{ex:subsec-gauge}

For a normalized weighted operator $\mathcal N$, let
$\mathcal B_{\mathcal N}(z)$ denote the self-adjoint realization of
the bipartization of $\mathcal N-zI$, and define its root resolvent
block by
\begin{equation}\label{ex:eq-gauge-root-compression}
  R^{\mathcal N}_{oo}(z,\eta)
  :=
  \Pi_o\bigl(\mathcal B_{\mathcal N}(z)-\eta I\bigr)^{-1}\Pi_o^*,
  \qquad \eta\in\mathbb C_+.
\end{equation}

The phase-removal and shift-rotation arguments are both consequences
of the following deterministic observation.

\begin{lemma}[Gauge covariance of the root resolvent]
\label{ex:lem-root-gauge}
Let $\mathcal N$ and $\widetilde{\mathcal N}$ be weighted operators on
the same rooted vertex set, and assume that their shifted
bipartizations have self-adjoint realizations with
$\mathcal D_{\mathrm{fin}}(V^b)$ as a common core.
Suppose that $U_+$ and $U_-$ are diagonal unitaries preserving
$\mathcal D_{\mathrm{fin}}(V)$ and that, on this domain,
\begin{equation}\label{ex:eq-general-gauge}
  \widetilde{\mathcal N}-\widetilde zI
  =
  U_+^*(\mathcal N-zI)U_-.
\end{equation}
If
  $U_+e_o=\gamma_+e_o$, $U_-e_o=\gamma_-e_o$ and $
  E_o=
  \begin{pmatrix}
    \gamma_+&0\\
    0&\gamma_-
  \end{pmatrix}$,
then
\begin{equation}\label{ex:eq-root-gauge-covariance}
  R^{\widetilde{\mathcal N}}_{oo}(\widetilde z,\eta)
  =
  E_o^*R^{\mathcal N}_{oo}(z,\eta)E_o,
  \qquad \eta\in\mathbb C_+.
\end{equation}
\end{lemma}

\begin{proof}
Put $\mathbb U=
  \begin{pmatrix}
    U_+&0\\
    0&U_-
  \end{pmatrix}$.
Equation~\eqref{ex:eq-general-gauge} and its adjoint give, on
$\mathcal D_{\mathrm{fin}}(V^b)$,
\[
  \mathcal B_{\widetilde{\mathcal N}}(\widetilde z)
  =
  \mathbb U^*\mathcal B_{\mathcal N}(z)\mathbb U.
\]
Since the finite-support domain is a core for both self-adjoint
operators, the identity extends to their self-adjoint realizations.
Consequently,
\[
  \bigl(
    \mathcal B_{\widetilde{\mathcal N}}(\widetilde z)-\eta I
  \bigr)^{-1}
  =
  \mathbb U^*
  \bigl(
    \mathcal B_{\mathcal N}(z)-\eta I
  \bigr)^{-1}
  \mathbb U.
\]
Moreover,
\[
  \mathbb U\Pi_o^*=\Pi_o^*E_o,
  \qquad
  \Pi_o\mathbb U^*=E_o^*\Pi_o.
\]
Compressing the preceding resolvent identity to the root pair proves
\eqref{ex:eq-root-gauge-covariance}.
\end{proof}

The next result states that the probability measure $\nu_{\alpha,\theta;z}$ depends only on
$\alpha$ and $|z|$. We henceforth denote it by
$\nu_{\alpha,z}^{\mathrm{SN}}$.
\begin{proposition}[Universality of the shifted singular value laws]
\label{ex:prop-gauge}
For every fixed $z\in\mathbb C$ and
$\eta\in\mathbb C_+$, the joint law of all four entries of the
normalized-PWIT root resolvent block is independent of the angular
measure $\theta$.
\end{proposition}

\begin{proof}
We need to address the phases and the shift. For simplicity, we suppress $\alpha$ from the notation.

\smallskip
\noindent\emph{Removal of the edge phases.}
Couple all angular laws using the same raw magnitudes and orientations.
For one realization, write $\mathcal M^\theta$ for the normalized
operator with its given edge phases, and let
$\mathcal M^{\mathrm{abs}}$ be obtained by replacing every nonzero raw
coefficient by its modulus.

Because the normalizers depend only on the edge magnitudes, the two
operators have the same normalizing factors. Since the underlying
graph is a tree and every unoriented edge carries exactly one nonzero
directed coefficient, we may recursively choose unit complex numbers
$(u_v)_{v\in V}$, with $u_o=1$, such that
\[
  \overline{u_u}\,m^\theta_{uv}\,u_v
  =
  |m^\theta_{uv}|
  \qquad\text{whenever }m^\theta_{uv}\ne0.
\]
Let $Ue_v=u_ve_v$. Then
\[
  U^*\mathcal M^\theta U=\mathcal M^{\mathrm{abs}}
\]
on finitely supported vectors. Applying
Lemma~\ref{ex:lem-root-gauge} with
  $U_+=U_-=U$ and $\widetilde z=z$,
gives the pathwise identity
\begin{equation}\label{ex:eq-gauge-root-equality}
  R^{\mathcal M^\theta}_{oo}(z,\eta)
  =
  R^{\mathcal M^{\mathrm{abs}}}_{oo}(z,\eta).
\end{equation}
The phase-free operator depends only on the raw magnitudes and
orientations, whose law depends on $\alpha$ but not on $\theta$.
Therefore the joint law of the whole root block is independent of
the angular measure.

\medskip
\noindent\emph{Rotation of the shift.}
Fix one realization of $\mathcal M$. Since the support graph is a
tree with one directed coefficient on each edge, there exists an
integer-valued function $h:V\to\mathbb Z$, with $h(o)=0$, such that
\[
  h(u)-h(v)=1
  \qquad\text{whenever }m_{uv}\ne0.
\]
For $\vartheta\in\mathbb R$, define
  $V_\vartheta e_v
  :=
  e^{i\vartheta h(v)}e_v$.
Then we have
  $V_\vartheta^*\mathcal M V_\vartheta
  =
  e^{-i\vartheta}\mathcal M$, or equivalently,
\[
  \mathcal M-e^{i\vartheta}zI
  =
  V_\vartheta^*(\mathcal M-zI)
  \bigl(e^{i\vartheta}V_\vartheta\bigr).
\]
Now we apply Lemma~\ref{ex:lem-root-gauge} with
\[
  U_+=V_\vartheta,
  \qquad
  U_-=e^{i\vartheta}V_\vartheta,
  \qquad
  \widetilde z=e^{i\vartheta}z.
\]
Since $h(o)=0$, the corresponding root matrix is
$E_\vartheta=
  \begin{pmatrix}
    1&0\\
    0&e^{i\vartheta}
  \end{pmatrix}$.
Hence
\begin{equation}\label{ex:eq-gauge-rotated-root-block}
  R^{\mathcal M}_{oo}(e^{i\vartheta}z,\eta)
  =
  E_\vartheta^*
  R^{\mathcal M}_{oo}(z,\eta)
  E_\vartheta
\end{equation}
and in particular,
 $ \operatorname{Tr}
  R^{\mathcal M}_{oo}(e^{i\vartheta}z,\eta)
  =
  \operatorname{Tr}
  R^{\mathcal M}_{oo}(z,\eta)$.
Finally, the root-resolvent characterization gives
\[
  \int_{\mathbb R}
  \frac1{x-\eta}\,
  \check\nu_{\alpha,\theta;z}(dx)
  =
  \frac12\mathbb E
  \operatorname{Tr}
  R^{\mathcal M^\theta}_{oo}(z,\eta).
\]
Equation~\eqref{ex:eq-gauge-root-equality} removes the dependence on
$\theta$, while
\eqref{ex:eq-gauge-rotated-root-block} removes the dependence on the
argument of $z$. Uniqueness of the Stieltjes transform, followed by
push-forward under $x\mapsto|x|$, proves that
$\nu_{\alpha,\theta;z}$ depends only on $\alpha$ and $|z|$.
\end{proof}

\begin{remark}[The raw PWIT]
Lemma~\ref{ex:lem-root-gauge} applies in exactly the same way with
the normalized operator $\mathcal M^\theta$ replaced by the raw
operator $\mathcal A^\theta$. Thus the raw root resolvent block is
also independent of the angular law after conditioning on the edge
magnitudes and orientations. This angular universality for the raw
shifted singular value laws is already established in
\cite[Theorem~1.1, Section~2.8, and Theorem~4.1]{BCC}.
For the normalized model, the additional observation is that the
normalizing factors depend only on magnitudes and therefore commute
with the phase-removal gauge.
\end{remark}
\subsection{Deterministic singular value inequalities}
\label{ex:subsec-deterministic}

For a matrix $M$ with rows $R_i$, write
$W_i:=\operatorname{span}_{\mathbb C}\{R_j:j\ne i\}$ and
$\Delta_i:=\operatorname{dist}(R_i,W_i)$.
All subspaces below are complex linear subspaces, also when the entries
are real.

\begin{lemma}[Row distances and negative moments]
\label{ex:lem-linear-algebra}
For an $n\times n$ matrix $M$, it holds
\begin{equation}\label{ex:eq-row-bound}
 \min_i\Delta_i\le\sqrt n\,s_n(M).
\end{equation}
For an $r\times n$ matrix $M$ with rank $r$, it holds for $0<p\le2$, 
\begin{equation}\label{ex:eq-fractional-negative}
 \sum_{j=1}^r s_j(M)^{-p}\le\sum_{i=1}^r\Delta_i^{-p}.
\end{equation}
If $M'$ is obtained from a square $M$ by deleting $k$ rows, then
$s_j(M)\ge s_j(M')\ge s_{j+k}(M)$ at every admissible index.
\end{lemma}

\begin{proof}
The inequality \eqref{ex:eq-row-bound} can be found in Lemma 4.16 of \cite{bordenave2012around}. The inequality \eqref{ex:eq-fractional-negative} for $p=2$ can be found in Lemma 4.14 of \cite{bordenave2012around}. For general $p \in (0,2]$  let
$G=MM^*$. The Schur complement gives $(G^{-1})_{ii}=\Delta_i^{-2}$,
while the eigenvalues of $G^{-1}$ are $s_j(M)^{-2}$.
Schur--Horn and concavity of $x\mapsto x^{p/2}$ show that the sum of
this function over the diagonal is at least its sum over the
eigenvalues. 
The last assertion follows by the Cauchy interlacing theorem for singular values.
\end{proof}

We also have, for $0<p\le2$, that
\begin{equation}\label{ex:eq-moments-deterministic}
 \int s^p\,\nu_{Y_n-zI}(ds)\le(1+|z|)^p
 \qquad \text{ and } \qquad
 \int|w|^p\,\mu_{Y_n}(dw)\le1.
\end{equation}
 For the eigenvalue
bound, apply Weyl's eigenvalue--singular value inequality at $z=0$;
see \cite[Lemma~B.5]{BCC}. Neither estimate requires centering or
symmetry. In particular all eigenvalue measures are deterministically
tight, since $\mu_{Y_n}\{|w|>R\}\le R^{-2}$.

\subsection{Polynomial invertibility for real and complex entries}
\label{ex:subsec-invertibility}
For the remainder of this section, we additionally impose the bounded density assumption \eqref{eq:model-density}.

The common denominator in \eqref{eq:model-row-normalization}  makes a self-normalized row of $X_n$ a nonlinear function of
its independent coordinates. We first give a deterministic one-variable
estimate, and then apply it either to a real density or to sections of
a joint planar density. 

\begin{lemma}[A real sublevel-length estimate]\label{ex:lem-sublevel}
There is a universal constant $C$ such that, for $a>0$, $b,c\in\mathbb R$,
$q>0$, $T\ge1$, and $u>0$,
\begin{equation}\label{ex:eq-sublevel}
 \begin{split}
 &\operatorname{Leb}\left\{x\in\mathbb R:|x|\le T\sqrt q,
       \left|\frac{ax+b}{\sqrt{x^2+q}}-c\right|\le u\right\}\\
 &\hspace{16mm}\le C\sqrt q(1+T^2)
                 \min\{1,\sqrt{u/a}\}.
 \end{split}
\end{equation}
\end{lemma}

\begin{proof}
Write $h(x)=(ax+b)/\sqrt{x^2+q}$ and
$A=(a^2+b^2/q)^{1/2}\ge a$. Direct calculation gives
\[
 h'(x)=\frac{aq-bx}{(x^2+q)^{3/2}},\qquad
 A^2-h(x)^2=\frac{(aq-bx)^2}{q(x^2+q)}.
\]
There are at most two monotonicity intervals, and on $|x|\le T\sqrt q$,
\[
 |h'(x)|=\frac{\sqrt q}{x^2+q}\sqrt{A^2-h(x)^2}
 \ge\frac{\sqrt{A^2-h(x)^2}}{(1+T^2)\sqrt q}.
\]
Changing variables separately on the monotonicity intervals (and taking
limits at their endpoints) bounds the length by
\[
 2(1+T^2)\sqrt q
       \int_{[c-u,c+u]\cap(-A,A)}\frac{ds}{\sqrt{A^2-s^2}}.
\]
The integral over the whole interval is $\pi$. If $u<A$, split it at
zero. On either half the largest integral over an interval of length
at most $2u$ occurs adjacent to the endpoint. When $2u\le A$,
\[
 \int_{A-2u}^A\frac{ds}{\sqrt{A^2-s^2}}
 \le A^{-1/2}\int_{A-2u}^A(A-s)^{-1/2}\,ds
       \le C\sqrt{u/A}.
\]
When $2u>A$, the same bound follows from $\pi/2$ after increasing $C$.
The integral is therefore at most $C\min(1,\sqrt{u/A})$.
Use $A\ge a$ to conclude.
\end{proof}
\begin{proposition}[Small balls for a normalized row]
\label{ex:prop-small-ball}
Let $U=(X_1,\ldots,X_n)/(\sum_{j=1}^n |X_j|^2)^{1/2}$, with $(X_j)_j$ independent
copies of the nondegenerate random variable $X$. For every $D_0>0$ there are constants $C_{D_0}$ and
$\kappa_{D_0}$ such that, for all sufficiently large $n$ and $0<u\le1$,
\begin{equation}\label{ex:eq-small-ball}
 \sup_{\|v\|_2=1}\sup_{w\in\mathbb C}
 \mathbb P\left(\left|\sum_{j=1}^n U_j\overline{v_j}-w\right|\le u\right)
 \le C_{D_0}n^{\kappa_{D_0}}u^{1/2}+C_{D_0}n^{-D_0}.
\end{equation}
In particular, for every $D_0>0$ some $b_{D_0}<\infty$ satisfies
\begin{equation}\label{ex:eq-small-ball-poly}
 \sup_{\|v\|_2=1,w\in\mathbb C}
 \mathbb P\left(\left|\sum_{j=1}^n U_j\overline{v_j}-w\right|\le n^{-b_{D_0}}\right)
 \le n^{-D_0}
\end{equation}
for all sufficiently large $n$. 
\end{proposition}

\begin{proof}
Choose $j$ with $a=|v_j|\ge n^{-1/2}$, and a unit complex number $\lambda$
so that $\lambda\overline{v_j}=a$. The event in
\eqref{ex:eq-small-ball} implies its real projection after multiplication
by $\lambda$. Conditional on the other coordinates, put
\[
 q=\sum_{k\ne j}|X_k|^2,\qquad
 b=\operatorname{Re}\left(\lambda\sum_{k\ne j}X_k\overline{v_k}\right),
 \qquad c=\operatorname{Re}(\lambda w).
\]
The projected event for real $X_j=x$ is
$|(ax+b)/\sqrt{x^2+q}-c|\le u$.
For complex $X_j=x+iy$, it is
$|(ax+b)/\sqrt{x^2+y^2+q}-c|\le u$.
Only the other coordinates have been conditioned on; $X_j$ retains its
original real or joint planar density.

Fix $0<\rho<\alpha$. Then regular variation \eqref{eq:model-radial-tail} and the Potter bounds give that 
$\mathbb P(|X|>s)\le  s^{-\rho}$ for sufficiently large $s>0$. Moreover since $X$ is non-degenerate we get that there is some $b_0>0$ with $\mathbb P(|X|\ge b_0)>0$.
A binomial lower-tail estimate therefore gives
$\mathbb P(q<c_0n)\le e^{-c_1n}$ for some $c_0,c_1>0$.
Choose $M,H>0$ so large that
\[
 1-\rho M\le-(D_0+3),\qquad
 \rho(H+1/2)>D_0+3,
\]
and put $T=n^H$. Then
\[
 \mathbb P(q>n^{2M+1})\le n\mathbb P(|X|>n^M)
       \le Cn^{-(D_0+3)},
\]
and on  $\{q\ge c_0n\}$,
\[
 \mathbb P(|X_j|>T\sqrt q\mid X_k,\ k\ne j)
       \le Cn^{-(D_0+3)}.
\]
It remains to estimate the event on
$c_0n\le q\le n^{2M+1}$ and $|X_j|\le T\sqrt q$.

For real entries, Lemma~\ref{ex:lem-sublevel} and the density bound
$K=\|f\|_\infty$ give conditional probability at most
\[
 CK\sqrt q(1+T^2)\sqrt{u/a}
       \le Cn^{M+2H+3/4}u^{1/2}.
\]
For complex entries, integrate the joint density over the planar
sublevel set. For fixed $y$ with $|y|\le T\sqrt q$, apply
Lemma~\ref{ex:lem-sublevel} with $q_y=q+y^2$.
The allowed $x$ satisfy $|x|\le T\sqrt q\le T\sqrt{q_y}$, so their
length is at most
$C\sqrt{q_y}(1+T^2)\sqrt{u/a}$.
Since $q_y\le q(1+T^2)$, integration in $y$ bounds the area by
\[
 CTq(1+T^2)^{3/2}\sqrt{u/a}.
\]
Multiplying by $K$ gives conditional probability at most
$Cn^{2M+4H+5/4}u^{1/2}$. This is a bound using the joint density, not
the conditional density of $x$ given $y$.

Combining the estimates proves \eqref{ex:eq-small-ball}, with for
example $\kappa_{D_0}=2M+4H+5/4$ for both cases. To obtain
\eqref{ex:eq-small-ball-poly}, first use \eqref{ex:eq-small-ball} with
failure exponent $D_0+2$, then choose $b_{D_0}$ large enough. The extra
power absorbs the constants and the residual term.
\end{proof}

\begin{proposition}[Polynomial least singular value]
\label{ex:prop-invertibility}
For every $D_0>0$ there exists $B=B(D_0)<\infty$ such that, for all
sufficiently large $n$,
\begin{equation}\label{ex:eq-invertibility}
 \sup_{z\in\mathbb C}\mathbb P(s_n(Y_n-zI)\le n^{-B})\le n^{-D_0}.
\end{equation}
%The supremum is over deterministic numerical values of $z$.
\end{proposition}

\begin{proof}
Let $R_i=U_i-ze_i$ be row $i$, and $W_i$ the span of the other rows of $Y_n-zI$.
Choose a measurable unit vector $v_i\in W_i^\perp$, depending only on
those other rows; measurable Gram--Schmidt with deterministic tie
choices is one way to do so. Therefore, $v_i$ is independent of the raw row producing
$U_i$. Then using \eqref{ex:eq-row-bound}
\[
 \min_i \left|\sum_j(U_i)_j\overline{(v_i)_j}-z\overline{(v_i)_i}\right|
       \le \min_i \operatorname{dist}(R_i,W_i)\leq \sqrt{n}s_n(Y_n-zI)\,.
\]
Condition on $v_i$ and apply \eqref{ex:eq-small-ball-poly}, with an
arbitrary center and failure exponent $D_0+2$. If its threshold is
$n^{-b}$, Lemma~\ref{ex:lem-linear-algebra} and a union bound give
\[
 \mathbb P(s_n(Y_n-zI)\le n^{-b-1/2})\le n^{-D_0-1}.
\]
This proves the result after renaming the exponent.
\end{proof}

\subsection{Regular variation and a high-codimension distance estimate}
\label{ex:subsec-distance}

The following replacement for the pure-power inverse distance input
retains the quantile $a_m$ at the codimension scale. This is essential
for allowing an arbitrary slowly varying tail. We use the truncation
and convex concentration method of \cite[Proposition~3.3 and
Theorem~3.4]{BCC}.
\begin{lemma}[Radial moment and truncation estimates]
\label{ex:lem-radial}
Let $\mathbf X=(X_1,\ldots,X_n)$ and $L_n=\|\mathbf X\|_2$.
For every $0<s<\alpha$ and all sufficiently large $n$,
\begin{equation}\label{ex:eq-row-norm-moment}
  \mathbb E L_n^s\le C_s a_n^s.
\end{equation}
For every $\eta>0$ and all sufficiently large $n\ge m$,
\begin{equation}\label{ex:eq-quantile-ratio}
  \frac{a_n}{a_m}
  \le C_\eta\left(\frac nm\right)^{1/\alpha+\eta}.
\end{equation}
If
\[
  \ell_2(b)
  =\mathbb E\bigl[|X|^2\mathbf1_{\{|X|\le b\}}\bigr],
  \qquad
  m_b=\mathbb E[X\mid |X|\le b],
\]
then $\ell_2(b)\to\infty$ and
\begin{equation}\label{ex:eq-truncated-variance}
  \tau_b^2
  :=\mathbb E\bigl[|X-m_b|^2\mid |X|\le b\bigr]
  \sim\ell_2(b),
  \qquad
  \tau_b^2\ge c b^2\mathbb P(|X|>b)
\end{equation}
for all sufficiently large $b$.
\end{lemma}

\begin{proof}
Write
\[
  \overline F(b)=\mathbb P(|X|>b).
\]
We use Potter's theorem, Karamata's integration theorem, and the
asymptotic inversion theorem in
\cite[Theorems~1.5.6, 1.5.11, and~1.5.12]{BGT87}.

Since $\overline F\in\mathrm{RV}_{-\alpha}$ and $0<\alpha<2$,
integration by parts and Karamata's theorem give
\begin{equation}\label{ex:eq-karamata-truncated-second}
\begin{aligned}
  \ell_2(b)
  &=2\int_0^b t\overline F(t)\,dt-b^2\overline F(b)\\
  &\sim \frac{\alpha}{2-\alpha}\,
          b^2\overline F(b)
  \longrightarrow\infty.
\end{aligned}
\end{equation}
Moreover, asymptotic inversion implies that $(a_n)$ is regularly
varying with index $1/\alpha$. Potter's theorem therefore yields
\eqref{ex:eq-quantile-ratio}.

For the row-norm estimate, put
\[
  Z_j=|X_j|^2,
  \qquad
  \overline F_Z(u)=\mathbb P(Z_1>u)
                  =\overline F(\sqrt u).
\]
Then $\overline F_Z\in\mathrm{RV}_{-\alpha/2}$, and
\eqref{ex:eq-karamata-truncated-second} gives
\[
  \mathbb E\bigl[Z_1\mathbf1_{\{Z_1\le u\}}\bigr]
  \le C u\overline F_Z(u)
\]
for all sufficiently large $u$. Hence, by splitting according to
whether $\max_{j\le n}Z_j>u$ and applying Markov's inequality to the
truncated sum,
\[
\begin{aligned}
  \mathbb P(L_n^2>u)
  &\le n\overline F_Z(u)
     +\frac{n}{u}
       \mathbb E\bigl[Z_1\mathbf1_{\{Z_1\le u\}}\bigr]\\
  &\le Cn\overline F_Z(u).
\end{aligned}
\]
Taking $u=a_n^2T^2$, using
$n\overline F(a_n)\to1$, and applying Potter's theorem, we obtain,
for every $\varepsilon\in(0,\alpha)$,
\[
  \mathbb P(L_n>a_nT)
  \le C_\varepsilon T^{-\alpha+\varepsilon},
  \qquad T\ge1,
\]
uniformly for all sufficiently large $n$.

If $0<s<\alpha$, choose $\varepsilon<\alpha-s$ and integrate the
tail:
\[
  \mathbb E(L_n/a_n)^s
  \le
  1+sC_\varepsilon
  \int_1^\infty T^{s-1-\alpha+\varepsilon}\,dT
  \le C_s.
\]
This proves \eqref{ex:eq-row-norm-moment}.

It remains to center the truncated law. Set
\[
  p_b=\mathbb P(|X|\le b),
  \qquad
  \mu_b=\mathbb E\bigl[X\mathbf1_{\{|X|\le b\}}\bigr],
\]
so that $m_b=\mu_b/p_b$. For every fixed $K<b$,
Cauchy--Schwarz gives
\[
  |\mu_b|
  \le
  \mathbb E\bigl[|X|\mathbf1_{\{|X|\le K\}}\bigr]
  +\ell_2(b)^{1/2}\mathbb P(|X|>K)^{1/2}.
\]
Divide by $\ell_2(b)^{1/2}$, let $b\to\infty$, and then let
$K\to\infty$. Since $\ell_2(b)\to\infty$, this proves
\[
  |\mu_b|^2=o(\ell_2(b)).
\]
As $p_b\to1$,
\[
  \tau_b^2
  =\frac{\ell_2(b)}{p_b}
   -\frac{|\mu_b|^2}{p_b^2}
  \sim\ell_2(b).
\]
The final lower bound in \eqref{ex:eq-truncated-variance} follows from
\eqref{ex:eq-karamata-truncated-second}.
\end{proof}
The following lemma is a slightly modified version of Proposition~3.7 of \cite{BCC}, which allows us to prove Theorem~\ref{thm:interior-convergence} under the stated assumptions rather than under Assumption~\textup{(H2)} of \cite{BCC}.
\begin{lemma}[Raw inverse distances at the quantile scale]
\label{ex:lem-raw-distance}
Fix $\gamma=1/16$. Let $V\subset\mathbb C^n$ have complex codimension
$m\ge n^{1-\gamma}$. There is an event $E_{n,V}$ such that, uniformly
over all such $m,V$ and all sufficiently large $n$,
\begin{equation}\label{ex:eq-raw-distance-failure}
 \mathbb P(E_{n,V}^c)\le C e^{-c n^\gamma},
\end{equation}
and, for every fixed $s>0$,
\begin{equation}\label{ex:eq-raw-distance-moment}
 \mathbb E[\operatorname{dist}(\mathbf X,V)^{-s};E_{n,V}]
        \le C_s a_m^{-s}.
\end{equation}
The event can be chosen measurably in $V$. 
\end{lemma}
\begin{proof}
We use the truncation, centering, and median-concentration argument in
\cite[proof of Proposition~3.3, equations~(3.4)--(3.8)]{BCC}, based on
its Theorem~3.4. We retain the truncation level as a free parameter and
record the estimates needed to pass to the codimension quantile.

Write $\overline F(b)=\mathbb P(|X|>b)$,
$D_V=\operatorname{dist}(\mathbf X,V)$, and
$N_b=\#\{j:|X_j|>b\}$. Similarly to the cited argument one may find constants $c,C,K>0$
such that, for all sufficiently large $b$ and $m\ge4$,
\begin{equation}\label{ex:eq-bcc-variable-cutoff}
 \mathbb P(D_V\le b)
 \le \mathbb P(N_b>m/4)
       +4\exp\!\left(-c\,m\tau_b^2/b^2\right),
 \qquad m\tau_b^2/b^2\ge K,
\end{equation}
where $\tau_b^2=\mathbb E[|X-m_b|^2\mid |X|\le b]$.

Recall $\gamma=1/16$. For $u_0\le u\le n^\gamma$, where $u_0$ is a
sufficiently large constant, choose
\begin{equation}\label{ex:eq-variable-truncation}
 b=a_m u^{-2/\alpha}.
\end{equation}
Since $m\ge n^{1-\gamma}$ and $a_m$ is regularly varying with index
$1/\alpha$, these levels tend to infinity uniformly in $m,u$.
Potter's bounds, with exponent slack $\alpha/2$, and
$m\overline F(a_m)\to1$ give
\begin{equation}\label{ex:eq-truncation-tail-window}
 c u\le m\overline F(b)\le C u^3.
\end{equation}
Lemma~\ref{ex:lem-radial} therefore implies
$m\tau_b^2/b^2\ge c u$. Moreover,
\[
 \mathbb EN_b=n\overline F(b)
       \le C(n/m)u^3\le Cn^{4\gamma}=o(m).
\]
The exponential Markov bound for this binomial variable yields,
uniformly over the same range,
\begin{equation}\label{ex:eq-deleted-coordinates}
 \mathbb P(N_b>m/4)
 \le\exp\{-m/4+(e-1)n\overline F(b)\}\le e^{-cm}.
\end{equation}
Taking $u_0$ large enough and applying
\eqref{ex:eq-bcc-variable-cutoff}, we conclude that
\begin{equation}\label{ex:eq-variable-small-distance}
 \mathbb P(D_V\le a_m u^{-2/\alpha})
       \le e^{-cm}+4e^{-cu},
       \qquad u_0\le u\le n^\gamma.
\end{equation}

Define
\begin{equation}\label{ex:eq-distance-event}
 E_{n,V}=\{D_V>a_m n^{-2\gamma/\alpha}\}.
\end{equation}
Taking $u=n^\gamma$ proves
\eqref{ex:eq-raw-distance-failure}. Set
$R=a_m/D_V$ on $E_{n,V}$ and $R=0$ otherwise. Then
$0\le R\le n^{2\gamma/\alpha}$, and
\eqref{ex:eq-variable-small-distance}, with $u=v^{\alpha/2}$, gives
\[
 \mathbb P(R>v)\le e^{-cm}+4e^{-cv^{\alpha/2}},
 \qquad u_0^{2/\alpha}\le v\le n^{2\gamma/\alpha}.
\]
For every fixed $s>0$, tail integration therefore gives
\[
 \mathbb ER^s
 \le u_0^{2s/\alpha}
   +4s\int_0^\infty v^{s-1}e^{-cv^{\alpha/2}}\,dv
   +n^{2\gamma s/\alpha}e^{-cm}
 \le C_s.
\]
This proves \eqref{ex:eq-raw-distance-moment}.
The event in \eqref{ex:eq-distance-event} is measurable in
$\mathbf X$ and the orthogonal projection onto $V$, so the asserted
measurability is preserved.
\end{proof}
\begin{proposition}[Fractional normalized inverse distances]
\label{ex:prop-normalized-distance}
Let $0<p<\alpha$, $W\subset\mathbb C^n$, and $i\le n$. Put
$V=W+\operatorname{span}_{\mathbb C}\{e_i\}$ and
$m=\operatorname{codim}V$. If $m\ge n^{1-\gamma}$, the event
$E_{n,V}$ of Lemma~\ref{ex:lem-raw-distance} satisfies, uniformly in
$z\in\mathbb C$,
\begin{equation}\label{ex:eq-normalized-distance}
 \mathbb E[\operatorname{dist}(U-ze_i,W)^{-p};E_{n,V}]
       \le C_p(a_n/a_m)^p.
\end{equation}
Its failure probability is bounded by $Ce^{-cn^\gamma}$.
Both assertions remain valid conditionally when $W$ is random and
independent of the raw row $\mathbf X$.
\end{proposition}

\begin{proof}
Since $V$ contains $e_i$ and is linear,
\[
 \operatorname{dist}(U-ze_i,W)
       \ge\operatorname{dist}(U,V)
       =D_V/L_n.
\]
Choose conjugate exponent $q,q'>1$ with $pq<\alpha$. 
H\"older, \eqref{ex:eq-row-norm-moment}, and
\eqref{ex:eq-raw-distance-moment} at exponent $pq'$ give
\[
 \mathbb E[L_n^pD_V^{-p};E_{n,V}]
 \le (\mathbb EL_n^{pq})^{1/q}
       (\mathbb E[D_V^{-pq'};E_{n,V}])^{1/q'}
 \le C_p a_n^p a_m^{-p}.
\]
Uniformity over deterministic subspaces and the measurable choice of
the event permit conditioning on an independent random $W$.
\end{proof}

\subsection{Intermediate singular values and logarithmic integrability}
\label{ex:subsec-log-ui}

We obtain logarithmic uniform integrability by the argument of
\cite[Section~3.4]{BCC}, using the distance and invertibility estimates
established above. We retain the resulting negative-moment estimate,
which will also be used in the analysis of the limiting measure.

\begin{proposition}[Uniform integrability of the logarithm]
\label{ex:prop-ui}
For every fixed deterministic $z\in\mathbb C$,
\eqref{ex:eq-main-ui} holds. Moreover, there exist $q>0$, $C<\infty$,
and $n_0$, independent of the numerical value of $z$, and events
$\mathcal G_n(z)$ such that
\begin{equation}\label{ex:eq-negative-moment}
 \mathbb P(\mathcal G_n(z)^c)\longrightarrow0,
 \qquad
 \sup_{n\ge n_0}
 \mathbb E\left[
   \int s^{-q}\,\nu_{Y_n-zI}(ds);
   \mathcal G_n(z)
 \right]\le C.
\end{equation}
The limiting measure in Proposition~\ref{ex:prop-weak-sv} satisfies
\begin{equation}\label{ex:eq-limit-negative}
 \int s^{-q}\,\nu_{\alpha,z}^{\mathrm{SN}}(ds)\le C,
 \qquad
 \int s^2\,\nu_{\alpha,z}^{\mathrm{SN}}(ds)\le(1+|z|)^2,
\end{equation}
and in particular,
$ \int|\log s|\,\nu_{\alpha,z}^{\mathrm{SN}}(ds)<\infty$.
Here $\mathbb E[Z;E]$ denotes integration over $E$, with the
integrand set to zero outside $E$.
\end{proposition}

\begin{proof}
We apply the row-deletion argument of \cite[Section~3.4]{BCC}, using
Propositions~\ref{ex:prop-normalized-distance}
and~\ref{ex:prop-invertibility} in place of the corresponding
distance and invertibility estimates there.
The inverse-second-moment step is replaced by the fractional
negative-moment inequality of Lemma~\ref{ex:lem-linear-algebra}.
Together with the quantile-ratio estimate in
Lemma~\ref{ex:lem-radial}, these inputs give
\eqref{ex:eq-negative-moment} for some sufficiently small $q>0$.
The conditioning on the other rows is valid because the normalized
rows are independent. The uniformity of the preceding estimates in
the deterministic shift gives the asserted uniformity of $q,C,n_0$.

Uniform integrability at zero follows from
\eqref{ex:eq-negative-moment}, Markov's inequality, and
$s^q|\log s|\to0$ as $s\downarrow0$.
At infinity, it follows from the deterministic bound
\[
 \int s^2\,\nu_{Y_n-zI}(ds)\le(1+|z|)^2.
\]
This proves \eqref{ex:eq-main-ui}.

To obtain \eqref{ex:eq-limit-negative}, apply
Proposition~\ref{ex:prop-weak-sv} to the bounded continuous functions
\[
 f_K(s)=\min\{s^{-q},K\},\qquad f_K(0)=K,
 \qquad
 g_K(s)=\min\{s^2,K\}.
\]
For fixed $K$, convergence in probability to the deterministic limit
also gives convergence of the expectations of these bounded integrals.
Moreover,
\[
 \mathbb E\int f_K\,d\nu_{Y_n-zI}
 \le C+K\mathbb P(\mathcal G_n(z)^c).
\]
Letting $n\to\infty$ and then $K\uparrow\infty$ proves the first
bound in \eqref{ex:eq-limit-negative}. The second follows similarly
from the deterministic second-moment bound.
These two moment bounds imply logarithmic integrability of the
limiting measure.
\end{proof}

\begin{proof}[Completion of the proof of Theorem~\ref{thm:interior-convergence}]
Part~\textup{(i)} follows from
Propositions~\ref{ex:prop-weak-sv} and~\ref{ex:prop-gauge}, and
part~\textup{(ii)} is Proposition~\ref{ex:prop-ui}.
For part~\textup{(iii)}, apply the Hermitization argument of
\cite[Lemma~A.2]{BCC} to the shifted singular value convergence and
logarithmic uniform integrability established in
parts~\textup{(i)} and~\textup{(ii)}.
The empirical eigenvalue measures are deterministically tight, since
\[
 \int_{\mathbb C}|w|^2\,\mu_{Y_n}(dw)
 \le \frac1n\|Y_n\|_F^2=1.
\]
This gives convergence in probability of $\mu_{Y_n}$ and identifies
its limit as the distributional Laplacian of the limiting logarithmic
potential, as asserted.
\end{proof}
\section{Proof of the geometric and endpoint properties}
\label{geo:sec-proof}

The aim of this section is to prove Theorem~\ref{thm:limit-properties} and
Theorem~\ref{thm:parameterendpoints}. Throughout,
\begin{align}\label{eq:defn:parameters}
 0<\alpha<2,\qquad \beta:=\alpha/2,\qquad
 \chi_\alpha:=1-\beta,\qquad
 F_\alpha(r):=\mu_\alpha^{\mathrm{SN}}\{w:|w|\le r\}.
\end{align}

Write $\mathcal M=\mathcal M_{\alpha,\theta}$ and
$\mathcal B(z)=\mathcal B_{\alpha,\theta}(z)$ for the normalized PWIT
operator and its self-adjoint bipartization, with root $o=\varnothing$.
The plus coordinates are the row coordinates of $\mathcal M$, and the
minus coordinates are its column coordinates. We retain the notation
$m_{uv}=a_{uv}/d(v)$ from the existence proof. Its two structural
properties are
\begin{equation}\label{geo:eq-basic-weights}
 \sum_u|m_{uv}|^2=1,\qquad
 \sum_u|m_{vu}|^2<\infty\quad\text{almost surely for every }v,
\end{equation}
and exactly one of $a_{uv},a_{vu}$ is nonzero on each raw tree edge.
These follow from \eqref{ex:eq-M-definition} and
Lemma~\ref{ex:lem-reverse}.

\subsection{Rotational invariance and the support inclusion}
\label{geo:subsec-support}

\begin{lemma}[Rotations and a symmetric representative]
\label{geo:lem-support}
The measure $\mu_\alpha^{\mathrm{SN}}$ is rotationally invariant and
\[
 \operatorname{supp}\mu_\alpha^{\mathrm{SN}}
 \subseteq\overline{\mathbb D}.
\]
\end{lemma}

\begin{proof}
 The invariance follows by Proposition \ref{ex:prop-gauge}, Theorem \ref{thm:interior-convergence} (iii), since the potential in \eqref{eq:main-logarithmic-potential} is radial, and
its distributional Laplacian is rotation invariant.

To prove the support inclusion, use universality in
Theorem~\ref{thm:interior-convergence} and choose the real symmetric
representative with density
\begin{equation}\label{geo:eq-pareto-representative}
 f_\alpha(x)=\frac\alpha2|x|^{-1-\alpha}
                         \mathbf1_{\{|x|\ge1\}}.
\end{equation}
This is a bounded density with the required tail index and angular
law $\tfrac12(\delta_{-1}+\delta_1)$. Its limiting eigenvalue measure
is the same $\mu_\alpha^{\mathrm{SN}}$ as for every distribution of $X$ allowed in
Theorem~\ref{thm:interior-convergence}. For its normalized matrices,
\cite[Corollary~2.6]{franccois2026characteristic} gives, for every
$\varepsilon>0$,
\begin{equation}\label{geo:eq-external-radius-upper}
 \mathbb P(\rho(Y_n)\ge1+\varepsilon)\longrightarrow0, \qquad \nto\,.
\end{equation}
Let $f\ge0$ be continuous and compactly supported in
$\mathbb C\setminus\overline{\mathbb D}$. Choose $\varepsilon>0$
so that its support lies in the set $\{|w|>1+\varepsilon\}$. Then
\[
 \mathbb E\int f\,d\mu_{Y_n}
 \le\|f\|_\infty\,
          \mathbb P(\rho(Y_n)\ge1+\varepsilon)\longrightarrow0.
\]
Weak convergence in probability to a deterministic measure, and
boundedness of the integrals, imply
$\int f\,d\mu_\alpha^{\mathrm{SN}}=0$.
The desired claim follows.
\end{proof}
\subsection{A result on supercritical fresh recursions}

In this subsection we establish a general criterion for the outer support
radius of a measure characterized by a fresh cavity recursion. We first
prove the supercritical instability estimate used for the lower bound.
No independence between the variables $Q_t$ and $D_t$ is required.

\begin{proposition}\label{prop:supercritical-fresh-recursion}
Let $(Q_t,D_t)_{t>0}$ be a family of random pairs such that
\[
    0<Q_t<\infty \quad \text{a.s.},
    \qquad
    0\leq D_t\leq 1 \quad \text{a.s.}
\]
Suppose that there exist an integer $N\geq 2$ and a random vector
    $A=(A_1,\ldots,A_N)\in(0,\infty)^N$
such that, for every $t>0$, one can construct i.i.d.\ copies
    $(Q_{1,t},D_{1,t}),\ldots,(Q_{N,t},D_{N,t})$
of $(Q_t,D_t)$, independent of $A$, together with a random variable
$\widehat Q_t\stackrel{d}{=}Q_t$, satisfying
\begin{equation}\label{eq:fresh-recursion-lower-bound}
    \widehat Q_t
    \geq
    \sum_{k=1}^N A_k Q_{k,t}(1-D_{k,t})
    \qquad \text{a.s.}
\end{equation}
Assume moreover that
\begin{equation}\label{eq:supercritical-smoothing-condition}
    \inf_{0\leq\theta\leq 1} m(\theta)>1,
    \qquad \text{ where }\quad
    m(\theta):=
    \mathbb E\sum_{k=1}^N A_k^\theta .
\end{equation}
Then
\begin{equation}\label{eq:supercritical-D-positive}
    \liminf_{t\downarrow 0}\mathbb E D_t>0.
\end{equation}
\end{proposition}

\begin{proof}
Suppose, to the contrary, that \eqref{eq:supercritical-D-positive} fails.
Since $D_t\geq 0$, there exists a sequence $t_n\downarrow0$ such that
\begin{equation}\label{eq:Dn-to-zero}
    \mathbb E D_{t_n}\longrightarrow0.
\end{equation}
For brevity, write $Q_n:=Q_{t_n}$ and $D_n:=D_{t_n}$.
For $b>0$, let
\[
    G_n(b):=\mathbb E e^{-Q_n/b}.
\]
Since $0<Q_n<\infty$ almost surely, the map $b\mapsto G_n(b)$ is continuous
and strictly increasing. Moreover, by dominated convergence, we have
\[
    \lim_{b\downarrow0}G_n(b)=0
    \qquad \text{ and } \qquad 
    \lim_{b\uparrow\infty}G_n(b)=1.
\]
Hence there exists a unique $b_n\in(0,\infty)$ such that
\begin{equation}\label{eq:laplace-normalization}
    \mathbb E e^{-Q_n/b_n}=\frac12.
\end{equation}
Set $Y_n:=\frac{Q_n}{b_n}$ and 
    $\phi_n(s):=\mathbb E e^{-sY_n}$ for $s\geq0$.
By construction,
\begin{equation}\label{eq:phi-normalization}
    \phi_n(1)=\frac12.
\end{equation}

We regard the law of $Y_n$ as a probability measure on the compact space
    $[0,\infty]:=[0,\infty)\cup\{\infty\}$, 
the one-point compactification of $[0,\infty)$.  Passing to a subsequence, the laws of $Y_n$ converge weakly
on $[0,\infty]$ to a probability measure $\nu$. With the convention
$e^{-s\infty}=0$, the function $x\mapsto e^{-sx}$ is continuous on
$[0,\infty]$ for every $s>0$, and therefore
\begin{equation}\label{eq:phi-limit}
    \phi_n(s)\longrightarrow
    \phi(s):=
    \int_{[0,\infty]} e^{-sx}\,\nu(dx),
    \qquad s>0.
\end{equation}
In particular, \eqref{eq:phi-normalization} gives
\begin{equation}\label{eq:phi-one-half}
    \phi(1)=\frac12.
\end{equation}

We next compare $Y_n$ with the damped variable that appears in the
recursion. Define
    $Z_n:=Y_n(1-D_n)$ and $\psi_n(s):=\mathbb E e^{-sZ_n}$.
For every $x\geq0$ and $d\in[0,1]$, convexity of
$d\mapsto e^{-x(1-d)}$ gives
\[
    e^{-x(1-d)}
    \leq
    (1-d)e^{-x}+d,
\]
and consequently
    $0\leq e^{-x(1-d)}-e^{-x}\leq d$.
It follows that
\begin{equation}\label{eq:psi-phi-close}
    \sup_{s\geq0}
    \left|\psi_n(s)-\phi_n(s)\right|
    \leq
    \mathbb E D_n
    \longrightarrow0.
\end{equation}
Notice that this estimate does not require any independence between
$Y_n$ and $D_n$.

Next, let
\[
    Y_{k,n}:=\frac{Q_{k,t_n}}{b_n} \qquad \text{ and } \qquad
    Z_{k,n}:=Y_{k,n}(1-D_{k,t_n}).
\]
Dividing \eqref{eq:fresh-recursion-lower-bound} by $b_n$, we obtain a random
variable $\widehat Y_n\stackrel{d}{=}Y_n$ such that
    $\widehat Y_n
    \geq
    \sum_{k=1}^N A_kZ_{k,n}$.
Since $x\mapsto e^{-sx}$ is decreasing,
\[
    \phi_n(s)
    \leq
    \mathbb E\exp\left\{
        -s\sum_{k=1}^N A_kZ_{k,n}
    \right\}.
\]
Conditioning on $A$, and using the independence of $A$ from the i.i.d.\
child pairs, gives
\begin{equation}\label{eq:phi-psi-recursion}
    \phi_n(s)
    \leq
    \mathbb E\prod_{k=1}^N\psi_n(sA_k).
\end{equation}
For numbers $a_k,b_k\in[0,1]$, we have the inequality
    $\left|\prod_{k=1}^N a_k-\prod_{k=1}^N b_k\right|
    \leq
    \sum_{k=1}^N|a_k-b_k|$.
Hence \eqref{eq:psi-phi-close} and \eqref{eq:phi-psi-recursion} imply
\[
    \phi_n(s)
    \leq
    \mathbb E\prod_{k=1}^N\phi_n(sA_k)
    +
    N\mathbb E D_n.
\]
Letting $n\to\infty$, using \eqref{eq:Dn-to-zero},
\eqref{eq:phi-limit}, and dominated convergence, we obtain
\begin{equation}\label{eq:smoothing-supersolution}
    \phi(s)
    \leq
    \mathbb E\prod_{k=1}^N\phi(sA_k),
    \qquad s>0.
\end{equation}

We claim that $\nu$ puts no mass at the point $\infty$. Letting
    $c:=\lim_{s\downarrow0}\phi(s)$,
monotone convergence yields
\[
    c=\nu([0,\infty))
      =1-\nu(\{\infty\}).
\]
Since $\phi$ is nonincreasing and $\phi(1)=1/2$, we have $c\geq1/2$.
Moreover, since $A_k>0$ almost surely, letting $s\downarrow0$ in
\eqref{eq:smoothing-supersolution} and applying dominated convergence yields
$c\leq c^N$.
As $N\geq2$ and $c>0$, this is possible only if $c=1$. Therefore
$\nu(\{\infty\})=0$. In particular, $\nu$ is a probability measure on
$[0,\infty)$ and $\phi$ is its Laplace transform. It is nonconstant, since
    $\phi(0)=1$ and $\phi(1)=\frac12$.

We now place \eqref{eq:smoothing-supersolution} in the standard framework
of the smoothing transform. Define
\[
    (\mathcal Sf)(s)
    :=
    \mathbb E\prod_{k=1}^N f(sA_k).
\]
Then \eqref{eq:smoothing-supersolution} reads
 $\phi\leq\mathcal S\phi$.
Since $\mathcal S$ is order preserving,
\begin{equation}\label{eq:smoothing-iterate}
    \phi\leq\mathcal S^j\phi,
    \qquad j\geq1.
\end{equation}

To express the iterates in branching notation, attach to each vertex $v$
of the rooted $N$-ary tree an independent copy
    $A(v)=(A_1(v),\ldots,A_N(v))$
of $A$, and define
\[
    L(\varnothing)=1,
    \qquad
    L(vk)=L(v)A_k(v).
\]
Thus $L(v)$ is the product of the edge weights along the path from the root
to $v$. By independence of the weight vectors attached to different
vertices,
\[
    (\mathcal S^j\phi)(s)
    =
    \mathbb E\prod_{|v|=j}\phi(sL(v)).
\]
Consequently, \eqref{eq:smoothing-iterate} becomes
\begin{equation}\label{eq:AM-branching-inequality}
    \phi(s)
    \leq
    \mathbb E\prod_{|v|=j}\phi(sL(v)),
    \qquad
    s\geq0,\quad j\geq1.
\end{equation}

The inequality \eqref{eq:AM-branching-inequality} is precisely the
inequality appearing in \cite[Lemma~6.3]{alsmeyer2012fixed}. To apply \cite[Lemma~6.3]{alsmeyer2012fixed}, define the weight sequence
\[
    T_k=A_k,\qquad 1\le k\le N,
    \qquad
    T_k=0,\qquad k>N.
\]
In the notation of \cite{alsmeyer2012fixed},
\[
    N_T:=\sum_{k\ge1}\mathbf 1_{\{T_k>0\}}=N
    \qquad\text{a.s.},
\]
and $m(\theta)
    =
    \mathbb E\sum_{k=1}^N A_k^\theta$.
Thus assumption (A2) of \cite{alsmeyer2012fixed} holds since
\(\mathbb E N_T=N>1\), while
\eqref{eq:supercritical-smoothing-condition} is precisely the strict
negation of condition (A3a) there, the latter being
\[
    \inf_{\theta\in[0,1]}m(\theta)\le1.
\]
Assumption (A1) of \cite{alsmeyer2012fixed} is
    $\mathbb P\bigl(T\in\{0,1\}^{\mathbb N}\bigr)<1$.
Since \(A_k>0\) almost surely, this condition can fail only if
\[
    A_1=\cdots=A_N=1
    \qquad\text{a.s.}
\]
In that exceptional case, however,
\eqref{eq:smoothing-supersolution} gives
$\phi(s)\le\phi(s)^N$,
which is impossible at \(s=1\), since
\(\phi(1)=1/2\) and \(N\ge2\). We may therefore assume that (A1)
holds.

Finally, \cite[Lemma~6.3]{alsmeyer2012fixed} requires 
\[
    \phi(s)
    \le
    \mathbb E\prod_{|v|=j}\phi(sL(v)),
    \qquad s\ge0,\ j\ge1,
\]
which is \eqref{eq:AM-branching-inequality}. Hence
\cite[Lemma~6.3]{alsmeyer2012fixed} implies
$\inf_{0\le\theta\le1}
    \mathbb E\sum_{k=1}^N A_k^\theta
    \le1$,
contradicting
\eqref{eq:supercritical-smoothing-condition}. Hence
\eqref{eq:Dn-to-zero} is impossible, and
\eqref{eq:supercritical-D-positive} follows.
\end{proof}

\begin{remark}\label{rem:infinite-fresh-recursion}
Proposition~\ref{prop:supercritical-fresh-recursion} applies immediately to an infinite fresh recursion. Suppose
that
\[
    Q_t
    =
    C_t+
    \sum_{k\geq1}
    A_k Q_{k,t}(1-D_{k,t}),
    \qquad
    C_t\geq0,
\]
where the child pairs are i.i.d.\ copies of $(Q_t,D_t)$ and are independent
of the entire weight sequence $(A_k)_{k\geq1}$. If, for some finite
$N\geq2$,
\[
    \inf_{0\leq\theta\leq1}
    \mathbb E\sum_{k=1}^N A_k^\theta>1,
\]
then \eqref{eq:fresh-recursion-lower-bound} follows by discarding the
nonnegative term $C_t$ and all terms with $k>N$. Hence
Proposition~\ref{prop:supercritical-fresh-recursion} yields
$\liminf_{t\downarrow0}\mathbb E D_t>0$.
\end{remark}
We continue with the following theorem.

\begin{theorem}[Outer support radius from a fresh cavity recursion]
\label{thm:general-support-criterion}
Let $d\in\{2,3,\ldots\}\cup\{\infty\}$, and let
$W=(W_k)_{1\le k\le d}$ be a random family such that
$0<W_k<\infty$ almost surely for every $k$. Set
    $c(p):=\mathbb E\sum_{k=1}^{d}W_k^p\in(0,\infty]$, $0<p\le1$,
and assume that $c(p_0)<\infty$ for some $p_0\in(0,1]$.
Define
\begin{equation}\label{eq:general-support-radius}
    \rho(W):=\inf_{0<p\le1}c(p)^{1/p}.
\end{equation}
Let $\mu$ be a rotationally invariant probability measure on $\mathbb C$,
and write
    $F(r):=\mu\{z:|z|\le r\}$
    as well as 
    $R:=\sup\{|z|:z\in\operatorname{supp}\mu\}$.
Suppose that, for every $r,t>0$, there is a random pair of finite
variables $P_{r,t},Q_{r,t}\ge t$. Define
\[
    H_{r,t}:=\frac{Q_{r,t}}{r^2+P_{r,t}Q_{r,t}},
    \qquad
    \widehat H_{r,t}:=\frac{P_{r,t}}{r^2+P_{r,t}Q_{r,t}},
    \qquad
    D_{r,t}:=\frac{P_{r,t}Q_{r,t}}{r^2+P_{r,t}Q_{r,t}},
\]
and assume the following properties.
\begin{enumerate}
\item For each $r,t>0$, if
$(P_{k,r,t},Q_{k,r,t})_{1\le k\le d}$ are independent copies of
$(P_{r,t},Q_{r,t})$, independent of the entire weight family $W$, then
\begin{equation}\label{eq:general-support-fresh}
    Q_{r,t}\stackrel{d}{=}
    t+\sum_{k=1}^{d}W_k
       \frac{Q_{k,r,t}}{r^2+P_{k,r,t}Q_{k,r,t}}.
\end{equation}
The law of $W$ is the same for all $r,t$.
\item For every $r,t>0$,
\begin{equation}\label{eq:general-support-trace}
    \mathbb E H_{r,t}=\mathbb E\widehat H_{r,t}.
\end{equation}
\item At every continuity point $r>0$ of $F$,
\begin{equation}\label{eq:general-support-recovery}
    1-F(r)=\lim_{t\downarrow0}\mathbb E D_{r,t}.
\end{equation}
\end{enumerate}
Then
\begin{equation}\label{eq:general-support-edge}
    R^2=\rho(W).
\end{equation}
More precisely,
\[
    \lim_{t\downarrow0}\mathbb E D_{r,t}=0
       \quad\text{if }r^2>\rho(W),
    \qquad
    \liminf_{t\downarrow0}\mathbb E D_{r,t}>0
       \quad\text{if }0<r^2<\rho(W).
\]
No independence among the weights $W_k$, or between $P_{r,t}$ and
$Q_{r,t}$, is required.
\end{theorem}

\begin{proof}
Fix $r>0$, suppress $r$ in the notation, and write
$P_t=P_{r,t}$, $Q_t=Q_{r,t}$, $H_t=H_{r,t}$,
$\widehat H_t=\widehat H_{r,t}$, and $D_t=D_{r,t}$.
The definitions imply
\begin{equation}\label{eq:general-support-basic}
    0<D_t<1,
    \qquad H_t,\widehat H_t\le t^{-1},
    \qquad H_t=\frac{Q_t}{r^2}(1-D_t),
    \qquad D_t=Q_t\widehat H_t.
\end{equation}

\emph{Upper bound.}
Let $p\in(0,1]$ satisfy $c(p)<r^{2p}$, and put
$\delta:=1-r^{-2p}c(p)>0$.
Subadditivity, the fresh recursion, and $H_t\le t^{-1}$ first give
\[
    \mathbb E Q_t^p\le t^p+t^{-p}c(p)<\infty.
\]
Thus the following moment inequality can be rearranged:
\[
    \mathbb E Q_t^p
    \le t^p+c(p)\mathbb E H_t^p
    \le t^p+r^{-2p}c(p)\mathbb E Q_t^p.
\]
Consequently,
\begin{equation}\label{eq:general-support-Q-moment}
    \mathbb E Q_t^p\le\frac{t^p}{\delta}.
\end{equation}
Writing $h_t:=\mathbb E H_t=\mathbb E\widehat H_t$, we also have
\begin{equation}\label{eq:general-support-H-moment}
    h_t\le t^{p-1}\mathbb E H_t^p
        \le r^{-2p}t^{p-1}\mathbb E Q_t^p
        \le\frac{r^{-2p}}{\delta}t^{2p-1}.
\end{equation}
It is not necessary that $h_t$ tend to zero.
Set $a:=p/(1+p)$. Since $0<D_t<1$, H\"older's inequality and
\eqref{eq:general-support-basic} yield
\begin{align*}
    \mathbb E D_t
    &\le\mathbb E D_t^a
      =\mathbb E\bigl[Q_t^a\widehat H_t^a\bigr]\\
    &\le(\mathbb E Q_t^p)^{1/(1+p)}
         (\mathbb E\widehat H_t)^{p/(1+p)}\\
    &\le\frac{r^{-2p^2/(1+p)}}{\delta}
              t^{2p^2/(1+p)}
      \longrightarrow0.
\end{align*}
This holds for every $p\in(0,1]$ with $c(p)<r^{2p}$.
For each $r^2>\rho(W)$, such an exponent exists by
\eqref{eq:general-support-radius}. Radial recovery therefore gives
$F(r)=1$ at every continuity radius above $\sqrt{\rho(W)}$.
Taking these radii down to $\sqrt{\rho(W)}$ proves
$R^2\le\rho(W)$.

\emph{Lower bound.}
Fix $0<r^2<\rho(W)$ and put $A_k:=W_k/r^2$.
For $0<\vartheta\le1$,
\[
    m(\vartheta):=\mathbb E\sum_{k=1}^{d}A_k^\vartheta
       =r^{-2\vartheta}c(\vartheta)>1,
\]
while $m(0)=d\ge2$, with the value $+\infty$ allowed.
For an integer $n\ge2$, set $d_n:=\min\{d,n\}$ and
\[
    A_k^{(n)}:=A_k\wedge n\quad(1\le k\le d_n),
    \qquad
    m_n(\vartheta):=
        \mathbb E\sum_{k=1}^{d_n}(A_k^{(n)})^\vartheta.
\]
Each $m_n$ is finite and continuous on $[0,1]$, and
$m_n(\vartheta)\uparrow m(\vartheta)>1$ for every
$\vartheta\in[0,1]$.
The increasing open sets
$\{\vartheta\in[0,1]:m_n(\vartheta)>1\}$ therefore cover
$[0,1]$. Compactness and continuity imply that, for some finite $n$,
\begin{equation}\label{eq:general-support-supercritical}
    \inf_{0\le\vartheta\le1}m_n(\vartheta)>1.
\end{equation}

Using \eqref{eq:general-support-basic}, the fresh recursion becomes
\[
    Q_t\stackrel{d}{=}
    t+\sum_{k=1}^{d}A_kQ_{k,t}(1-D_{k,t}),
\]
where the pairs $(Q_{k,t},D_{k,t})$ are independent copies of
$(Q_t,D_t)$, independent of $W$.

Remark~\ref{rem:infinite-fresh-recursion} consequently gives
\[
    \liminf_{t\downarrow0}\mathbb E D_{r,t}>0.
\]
At continuity radii below $\sqrt{\rho(W)}$, radial recovery yields
$F(r)<1$. For an arbitrary $r<\sqrt{\rho(W)}$, choose a continuity
radius $r'$ with $r<r'<\sqrt{\rho(W)}$ and use
$1-F(r)\ge1-F(r')>0$. This proves $R^2\ge\rho(W)$ and completes
the proof.
\end{proof}

\begin{corollary}[Normalized fresh weights]
\label{cor:normalized-support-criterion}
Under the hypotheses of Theorem~\ref{thm:general-support-criterion},
suppose additionally that
    $\sum_{k=1}^{d}W_k=1$ almost surely.
Then $R=1$.
\end{corollary}

\begin{proof}
For $0<p\le1$, every weight belongs to $(0,1]$, and therefore
\[
    \sum_{k=1}^{d}W_k^p\ge\sum_{k=1}^{d}W_k=1.
\]
Thus $c(p)^{1/p}\ge1$, whereas $c(1)=1$.
It follows that $\rho(W)=1$.
\end{proof}
\subsection{The root block and the independent source average}
\label{geo:subsec-cavity}

Fix a real shift $r\ge0$ and a regularization $t>0$.
For a child $u$ of $v$, let $V_{u\to v}$ be the descendant subtree
rooted at $u$, obtained by deleting the underlying edge $(v,u)$.
Let $\mathcal B_{u\to v}(r)$ be the self-adjoint realization of the
restriction of the shifted bipartization to
$V_{u\to v}^{+}\sqcup V_{u\to v}^{-}$.
All coefficients are inherited from the full normalized operator:
the matching edges of weight $-r$ are retained, and deleting the
parent edge does not recompute any normalizer $d(x)$.
The self-adjointness and finite-support core follow from the argument
of Proposition~\ref{ex:prop-esa}: restriction preserves row-square
summability and leaves every source-column squared norm at most one.

Recall the notation from Subsection \ref{ex:subsec-resolvents}. Define the cavity root block by
\[
 G_{u\to v}(r,it)
 :=\Pi_u\bigl(\mathcal B_{u\to v}(r)-itI\bigr)^{-1}\Pi_u^*.
\]

By \cite[Lemma~2.2]{BCC}, its diagonal entries are purely imaginary
with positive imaginary parts bounded by $t^{-1}$. Thus we write
\begin{equation}\label{geo:eq-cavity-H}
 [G_{u\to v}(r,it)]_{11}=iH^+_{u\to v},\qquad
 [G_{u\to v}(r,it)]_{22}=iH^-_{u\to v},\qquad
 0<H^\pm_{u\to v}\le t^{-1}.
\end{equation}
The following is the specialization of the tree resolvent identity
\cite[Lemma~2.5]{BCC} to the normalized directed PWIT.

\begin{lemma}[Root self-energies]
\label{geo:lem-root-block}
Set
\begin{equation}\label{geo:eq-PQ}
 P=P_o(r,t)=t+\sum_{u\sim o}|m_{ou}|^2H^-_{u\to o},\qquad
 Q=Q_o(r,t)=t+\sum_{u\sim o}|m_{uo}|^2H^+_{u\to o}.
\end{equation}
The sums are almost surely finite, and
\begin{equation}\label{geo:eq-root-block}
 G_o(r,it):=\Pi_o\bigl(\mathcal B(r)-itI\bigr)^{-1}\Pi_o^*
 =\frac1{r^2+PQ}
   \begin{pmatrix}iQ&-r\\-r&iP\end{pmatrix}.
\end{equation}
In particular, with
$H_\pm=\operatorname{Im}[G_o(r,it)]_{\pm\pm}$,
\begin{equation}\label{geo:eq-root-H}
 H_+=\frac{Q}{r^2+PQ},\qquad
 H_-=\frac{P}{r^2+PQ},\qquad
 0<H_\pm\le t^{-1},
\end{equation}
and
\begin{equation}\label{geo:eq-Q-bounds}
 t\le Q\le t+t^{-1},\qquad
 t\le P\le t+t^{-1}\sum_u|m_{ou}|^2.
\end{equation}
\end{lemma}

\begin{proof}
Finiteness of the sums and \eqref{geo:eq-Q-bounds} follow from
\eqref{geo:eq-basic-weights} and \eqref{geo:eq-cavity-H}.
We apply \cite[Lemma~2.5]{BCC} pathwise to the normalized coefficients.
With
\[
 U(r,it)=\begin{pmatrix}it&r\\r&it\end{pmatrix},\qquad
 C_{ou}=\begin{pmatrix}0&m_{ou}\\\overline{m_{uo}}&0\end{pmatrix},
\]
that identity gives
\[
 G_o(r,it)
 =-\left(
 U(r,it)+\sum_{u\sim o}
 C_{ou}G_{u\to o}(r,it)C_{ou}^*
 \right)^{-1}.
\]
The series converges absolutely in matrix norm by square summability
of the root rows and columns and the bound
$\|G_{u\to o}(r,it)\|\le t^{-1}$.
Since exactly one of $m_{ou},m_{uo}$ is nonzero on each underlying
edge, \eqref{geo:eq-cavity-H} yields
\[
 C_{ou}G_{u\to o}(r,it)C_{ou}^*
 =i\begin{pmatrix}
 |m_{ou}|^2H^-_{u\to o}&0\\
 0&|m_{uo}|^2H^+_{u\to o}
 \end{pmatrix}.
\]
Consequently,
\[
 G_o(r,it)
 =-\begin{pmatrix}iP&r\\r&iQ\end{pmatrix}^{-1}
 =\frac1{r^2+PQ}\begin{pmatrix}iQ&-r\\-r&iP\end{pmatrix}.
\]
This proves \eqref{geo:eq-root-block} and the formulas in
\eqref{geo:eq-root-H}. Finally, $P,Q\ge t$ imply
$H_+\le P^{-1}\le t^{-1}$ and $H_-\le Q^{-1}\le t^{-1}$.
The argument is deterministic conditional on the weighted tree;
no independence or equality in law of the cavity entries is used.
\end{proof}
We next record the Poisson identities needed for the source-side
recursion and the boundary estimate. Each orientation-thinned cloud
of raw squared magnitudes has intensity
$\beta x^{-1-\beta}\,dx$, as in \eqref{ex:eq-squared-intensity}.

\begin{lemma}[Stable totals and normalized weights]
\label{geo:lem-stable}
Let $(J_k)$ be a Poisson cloud of intensity
$\beta x^{-1-\beta}\,dx$, put $S=\sum_kJ_k$, and set $\pi_k=J_k/S$.
With $K_\beta=\Gamma(1-\beta)$, one has $0<S<\infty$ almost surely and
\begin{align}
 \mathbb E e^{-uS}&=e^{-K_\beta u^\beta}, &&u\ge0,
                         \label{geo:eq-stable-LT}\\
 \mathbb ES^{-s}&=K_\beta^{-s/\beta}
       \frac{\Gamma(1+s/\beta)}{\Gamma(1+s)}, &&s>0,
                         \label{geo:eq-stable-negative}\\
 \mathbb ES^s&=K_\beta^{s/\beta}
       \frac{\Gamma(1-s/\beta)}{\Gamma(1-s)}, &&0<s<\beta.
                         \label{geo:eq-stable-positive}
\end{align}
For $p>\beta$,
\begin{equation}\label{geo:eq-weight-moments}
 \mathbb E\sum_k\pi_k^p
  =\frac{\Gamma(p-\beta)}{\Gamma(p)\Gamma(1-\beta)}.
\end{equation}
In particular,
\begin{equation}\label{geo:eq-chi}
 \sum_k\pi_k=1,\qquad
 \mathbb E\sum_k\pi_k^2=1-\beta=\chi_\alpha.
\end{equation}
\end{lemma}

\begin{proof}
The positivity and finiteness of $S$, together with
\eqref{geo:eq-stable-LT}, follow from the positive stable Poisson-sum
representation in \cite[Lemma~4.3]{BCC}.
The moment identities \eqref{geo:eq-stable-negative} and
\eqref{geo:eq-stable-positive} are the specializations of
\cite[equation~(160)]{PY97} to stable index $\beta$ and L\'evy-tail
constant $C=1$.

By \cite[Proposition~6]{PY97}, the decreasing rearrangement of
$(\pi_k)$ has the Poisson--Dirichlet distribution
$\operatorname{PD}(\beta,0)$.
Applying \cite[equation~(6)]{PY97} with these parameters and
$f(x)=x^p$ gives, for $p>\beta$,
\[
 \mathbb E\sum_k\pi_k^p
 =\frac{1}{\Gamma(\beta)\Gamma(1-\beta)}
   \int_0^1 x^{p-\beta-1}(1-x)^{\beta-1}\,dx
 =\frac{\Gamma(p-\beta)}{\Gamma(p)\Gamma(1-\beta)}.
\]
This proves \eqref{geo:eq-weight-moments};
\eqref{geo:eq-chi} follows from $\sum_k\pi_k=1$ and the case $p=2$.
\end{proof}

\begin{lemma}[Independent source averaging and equal expected diagonals]
\label{geo:lem-source-trace}
At the fresh root,
\begin{equation}\label{geo:eq-Q-rde}
 Q\stackrel d=t+\sum_{k\ge1}\pi_kH_{+,k},
\end{equation}
where $(H_{+,k})$ are i.i.d. copies of $H_+$, independent of the whole
normalized weight sequence $(\pi_k)$ in Lemma~\ref{geo:lem-stable}.
Moreover,
\begin{equation}\label{geo:eq-trace-identity}
 \mathbb EH_+=\mathbb EH_-.
\end{equation}
\end{lemma}
\begin{proof}
The identity for $Q$ in \eqref{geo:eq-PQ} has a simple independent-copy
interpretation. At the fresh root, the edges contributing to $Q-t$
satisfy $a_{uo}\ne0$, so their connecting coefficients belong to
column $o$, rather than to the corresponding child columns.
Their raw squared magnitudes form a source Poisson cloud $(J_k)$,
and
\[
 d(o)^2=\sum_kJ_k,\qquad
 |m_{u_ko}|^2=\frac{J_k}{\sum_\ell J_\ell}=\pi_k.
\]
Removing such a connecting edge leaves a fresh normalized subtree
at the child: no normalizer inside that subtree contains the
connecting coefficient. The descendant subtrees are independent
of one another and of the entire root source cloud. Consequently,
\begin{equation}
 Q\stackrel d=t+\sum_{k\ge1}\pi_kH_{+,k},
\end{equation}
where the $H_{+,k}$ are i.i.d. copies of $H_+$, independent of the
whole normalized weight sequence $(\pi_k)$ from
Lemma~\ref{geo:lem-stable}.
Indexing the source edges by decreasing magnitude preserves these
independence properties. The series converges absolutely, since
$H_{+,k}\le t^{-1}$ and $\sum_k\pi_k=1$.

Concerning \eqref{geo:eq-trace-identity}, let $M_n=Y_n^{\mathsf T}$, put $C_n=M_n-rI$, and write
$R_n=(B_n(r)-itI)^{-1}$.
Its two diagonal $n\times n$ blocks are
\[
 [R_n]_{++}=it(C_nC_n^*+t^2I)^{-1},
 \qquad
 [R_n]_{--}=it(C_n^*C_n+t^2I)^{-1}.
\]
Since $C_n$ is square, $C_nC_n^*$ and $C_n^*C_n$ have the same
eigenvalues, including multiplicities. Hence, for every realization,
\[
 \operatorname{Tr}[R_n]_{++}
 =
 \operatorname{Tr}[R_n]_{--}.
\]

Let $G_{n,j}(r,it)$ be the $2\times2$ resolvent block at vertex $j$.
Simultaneous permutation invariance of the law of $M_n$ gives,
for either sign,
\[
 \mathbb E\,\operatorname{Im}[G_{n,1}(r,it)]_{\pm\pm}
 =
 \frac1n\mathbb E\,\operatorname{Im}
                   \operatorname{Tr}[R_n]_{\pm\pm}.
\]
The two expectations on the left are therefore equal.
By Proposition~\ref{ex:prop-root},
$G_{n,1}(r,it)\Rightarrow G_o(r,it)$.
The uniform resolvent bound $t^{-1}$ allows passage to expectations,
and \eqref{geo:eq-trace-identity} follows. 
\end{proof}

The identity \eqref{geo:eq-trace-identity} concerns expectations only;
it does not assert equality in distribution of $H_+$ and $H_-$.

\subsection{Radial recovery from the regularized potential}
\label{geo:subsec-radial-recovery}

For $t>0$ define
\begin{equation}\label{geo:eq-potential-t}
 \Phi_{\alpha,t}(z):=\frac12\int_0^\infty
                \log(s^2+t^2)\,\nu_{\alpha,z}^{\mathrm{SN}}(ds), \ \ \ z\in \C
 \qquad \phi_{\alpha,t}(r):=\Phi_{\alpha,t}(r), \ \ \ r\in \R.
\end{equation}
Write $\phi_\alpha(r):=\Phi_\alpha(r)$ for the unregularized potential. 

\begin{lemma}[The radial distribution from the root block]
\label{geo:lem-radial-recovery}
For $t>0$, $\Phi_{\alpha,t}$ is radial and subharmonic, and
$\phi_{\alpha,t}$ is continuously differentiable on $(0,\infty)$.
For every $r>0$,
\begin{equation}\label{geo:eq-Ft}
  F_{\alpha,t}(r)
  :=r\phi_{\alpha,t}'(r)
  =
  \mathbb E
  \frac{r^2}{r^2+P_o(r,t)Q_o(r,t)}
  \in[0,1].
\end{equation}
At every continuity point $r>0$ of $F_\alpha(r)
  :=
  \mu_\alpha^{\mathrm{SN}}\{z:|z|\le r\}$,
one has
\begin{equation}\label{geo:eq-radial-recovery}
  F_\alpha(r)
  =
  \lim_{t\downarrow0}F_{\alpha,t}(r).
\end{equation}
\end{lemma}
\begin{proof}
Set
\[
  b_{\alpha,t}(z)
  :=
  \mathbb E[G_o(z,it)]_{12}.
\]
We use the resolvent-to-measure argument of
\cite[Theorem~2.1 and proof of Lemma~4.5]{BCC}.
For every fixed $t>0$, the finite-dimensional Jacobi identity,
root-block convergence, and the resolvent bound give
\begin{equation}\label{geo:eq-regularized-resolvent-measure}
  \mu_{\alpha,t}
  :=
  \frac1{2\pi}\Delta\Phi_{\alpha,t}
  =
  -\frac1\pi\partial_z b_{\alpha,t}
  \qquad\text{in }\mathcal D'(\mathbb C).
\end{equation}
The measure $\mu_{\alpha,t}$ is radial by
Proposition~\ref{ex:prop-gauge}.

By Proposition~\ref{ex:prop-ui}, for some $q>0$,
\[
  0\le
  \Phi_{\alpha,t}(z)-\Phi_\alpha(z)
  \le
  C_qt^q
  \int s^{-q}\nu_{\alpha,z}^{\mathrm{SN}}(ds)
  \le Ct^q.
\]
Hence
\[
  \mu_{\alpha,t}
  \Rightarrow
  \frac1{2\pi}\Delta\Phi_\alpha
  =
  \mu_\alpha^{\mathrm{SN}}
  \qquad(t\downarrow0).
\]

For a real shift $r>0$, Lemma~\ref{geo:lem-root-block} gives
\[
  b_{\alpha,t}(r)
  =
  -\mathbb E
  \frac{r}{r^2+P_o(r,t)Q_o(r,t)}.
\]
The radial form of
\eqref{geo:eq-regularized-resolvent-measure} therefore yields
\[
  \mu_{\alpha,t}\{|z|\le r\}
  =
  -r b_{\alpha,t}(r)
  =
  \mathbb E
  \frac{r^2}{r^2+P_o(r,t)Q_o(r,t)}.
\]
The first equality also shows that
$r\phi_{\alpha,t}'(r)=\mu_{\alpha,t}\{|z|\le r\}$.
Letting $t\downarrow0$ and using weak convergence at every continuity
radius of $F_\alpha$ proves
\eqref{geo:eq-Ft} and \eqref{geo:eq-radial-recovery}.
\end{proof}

\subsection{Moment identities and the outer support radius}
\label{geo:subsec-outer}

For the rest of the cavity argument, $r>0,t>0$ and $\alpha$ are fixed
unless explicitly stated otherwise. As in the beginning of subsection \ref{geo:subsec-cavity}, we write
\[
 P=P_o(r,t),\quad Q=Q_o(r,t),\quad
 q=\mathbb EQ,\quad v=\operatorname{Var}(Q),\quad
 D_{r,t}=\frac{PQ}{r^2+PQ}.
\]
In particular $0<D_{r,t}<1$, and radial
recovery gives
\begin{equation}\label{geo:eq-tail-recovery}
 1-F_\alpha(r)=\lim_{t\downarrow0}\mathbb ED_{r,t}
\end{equation}
at continuity radii, see Lemma \ref{geo:lem-radial-recovery}.

\begin{lemma}[First and second self-energy moments]
\label{geo:lem-moments}
The following identities hold:
\begin{align}
 q&=t+\mathbb EH_+, &&
       \label{geo:eq-q-mean}\\
 v&=\chi_\alpha\operatorname{Var}(H_+), &&
       \label{geo:eq-q-variance}\\
 \mathbb E[QD_{r,t}]&=(1-r^2)q+r^2t. &&
       \label{geo:eq-QD}
\end{align}
If $r^4>\chi_\alpha$, then
\begin{equation}\label{geo:eq-relative-variance}
 \frac{v}{q^2}\le\frac{\chi_\alpha}{r^4-\chi_\alpha}\qquad \text{ and } \qquad
 \frac{\mathbb EQ^2}{q^2}
       \le\frac1{1-\chi_\alpha/r^4}.
\end{equation}
\end{lemma}

\begin{proof}
All the moments in question are finite at fixed $t$, by
\eqref{geo:eq-Q-bounds}. Taking expectation in
\eqref{geo:eq-Q-rde} gives \eqref{geo:eq-q-mean}.
Conditionally on $(\pi_k)$, the mean of $Q$ is the constant
$t+\mathbb EH_+$ and its variance is
$\operatorname{Var}(H_+)\sum_k\pi_k^2$. The infinite-series formula
follows by bounded convergence from partial sums. Taking expectation
and using \eqref{geo:eq-chi} gives \eqref{geo:eq-q-variance}.
The pathwise relation $QD_{r,t}=Q-r^2H_+$ then proves
\eqref{geo:eq-QD}. Finally $H_+\le Q/r^2$, so
\[
 v\le\chi_\alpha\mathbb EH_+^2
      \le\frac{\chi_\alpha}{r^4}(q^2+v).
\]
Rearranging yields \eqref{geo:eq-relative-variance}.
\end{proof}

\begin{proposition}[Outer support radius; Theorem~\ref{thm:limit-properties}(ii)]
\label{geo:prop-outer}
Let $0<\alpha<2$. The outer support radius of
$\mu_\alpha^{\mathrm{SN}}$ is one. Moreover, for every $0<r<1$,
\begin{equation}\label{geo:eq-noncollapse}
    \liminf_{t\downarrow0}\mathbb E D_{r,t}>0 \qquad \text{ and } \qquad 
    1-F_\alpha(r)>0.
\end{equation}
Writing $q(r,t):=\mathbb E Q_o(r,t)$, one also has, for every fixed $0<r<1$,
\begin{equation}\label{geo:eq-t-over-q}
    \liminf_{t\downarrow0}q(r,t)>0 \qquad \text{ and } \qquad 
    \frac{t}{q(r,t)}\longrightarrow0,
    \quad t\downarrow0.
\end{equation}
\end{proposition}

\begin{proof}
We apply Theorem~\ref{thm:general-support-criterion} to the
fresh weights $W_k=\pi_k$. These weights are strictly positive,
their joint law is independent of $r,t$, and
$\sum_k\pi_k=1$ almost surely. For every $0<p\le1$,
\[
    \mathbb E\sum_{k\ge1}\pi_k^p
    \ge
    \mathbb E\sum_{k\ge1}\pi_k
    =1,
\]
with equality at $p=1$. Consequently,
\[
    \rho(\pi)
    =
    \inf_{0<p\le1}
    \left(\mathbb E\sum_{k\ge1}\pi_k^p\right)^{1/p}
    =1.
\]

The remaining hypotheses of the criterion hold with
\[
    P_{r,t}=P_o(r,t),
    \qquad
    Q_{r,t}=Q_o(r,t),
    \qquad
    H_{r,t}=H_+,
    \qquad
    \widehat H_{r,t}=H_-.
\]
Indeed, the root self-energy formulas give
$P_o(r,t),Q_o(r,t)\ge t$ and
\[
    H_+=\frac{Q_o(r,t)}{r^2+P_o(r,t)Q_o(r,t)},
    \qquad
    H_-=\frac{P_o(r,t)}{r^2+P_o(r,t)Q_o(r,t)}.
\]
Lemma~\ref{geo:lem-source-trace} supplies the fresh recursion
and the identity $\mathbb EH_+=\mathbb EH_-$, while
Lemma~\ref{geo:lem-radial-recovery} supplies radial recovery.
Since $\mu_\alpha^{\mathrm{SN}}$ is rotationally invariant,
Theorem~\ref{thm:general-support-criterion} therefore gives
the outer support radius one and the first part of
\eqref{geo:eq-noncollapse}.

To obtain the second part of \eqref{geo:eq-noncollapse} at an arbitrary
$0<r<1$, choose a continuity radius $r'\in(r,1)$ of
$F_\alpha$. By \eqref{geo:eq-tail-recovery} and
\eqref{geo:eq-noncollapse},
\[
    1-F_\alpha(r)
    \ge
    1-F_\alpha(r')
    =
    \lim_{t\downarrow0}\mathbb E D_{r',t}
    >0.
\]

Finally, fix $0<r<1$ and abbreviate
$Q=Q_o(r,t)$ and $q=q(r,t)$.
The bound $Q\le t+t^{-1}$ ensures that $q<\infty$.
Taking expectations in the fresh recursion and using
Lemma~\ref{geo:lem-source-trace}, we obtain
\[
    q=t+\mathbb EH_+=t+\mathbb EH_-.
\]
Since $D_{r,t}=QH_-$ and $0\le D_{r,t}\le1$,
Cauchy--Schwarz gives
\begin{equation}\label{geo:eq-q-lower-mechanism}
    \mathbb E D_{r,t}
    \le
    \mathbb E\sqrt{QH_-}
    \le
    \sqrt{\mathbb EQ\,\mathbb EH_-}
    =
    \sqrt{q(q-t)}
    \le q.
\end{equation}
Together with \eqref{geo:eq-noncollapse}, this implies
\[
    \liminf_{t\downarrow0}q(r,t)
    \ge
    \liminf_{t\downarrow0}\mathbb E D_{r,t}
    >0.
\]
Hence $t/q(r,t)\to0$, proving \eqref{geo:eq-t-over-q}.
\end{proof}
\subsection{Absence of point atoms}
\label{geo:subsec-atoms}

\begin{proposition}[Theorem \ref{thm:limit-properties} (iii)]
\label{geo:prop-atoms}
For each $\alpha\in(0,2)$ there are
$\kappa_\alpha>0$ and $C_\alpha<\infty$ such that
\begin{equation}\label{geo:eq-negative-radial}
 \int_{\mathbb C}|w|^{-\kappa_\alpha}
                    \,\mu_\alpha^{\mathrm{SN}}(dw)\le C_\alpha.
\end{equation}
Thus $F_\alpha(r)\le C_\alpha r^{\kappa_\alpha}$ for $r>0$, and
$\mu_\alpha^{\mathrm{SN}}$ has no point atoms.
\end{proposition}

\begin{proof}
Use the symmetric representative \eqref{geo:eq-pareto-representative}
if desired, so that constants can be denoted solely by their
$\alpha$-dependence. The estimate \eqref{ex:eq-negative-moment} at
$z=0$ supplies an exponent $\kappa_\alpha>0$ and events
$\mathcal G_n$ with probability tending to one such that
\begin{equation}\label{geo:eq-finite-negative-input}
 \sup_{n\ge n_0}\mathbb E\left[
 \mathbf1_{\mathcal G_n}\frac1n
       \sum_{j=1}^n s_j(Y_n)^{-\kappa_\alpha}\right]\le C_\alpha.
\end{equation}
The matrices are invertible on these good events. Apply the
eigenvalue--singular value moment inequality to $Y_n^{-1}$, as in
\cite[Lemma~B.5]{BCC}, to obtain
\[
 \frac1n\sum_j|\lambda_j(Y_n)|^{-\kappa_\alpha}
 \le\frac1n\sum_j s_j(Y_n)^{-\kappa_\alpha}.
\]
For $K>0$ define the bounded continuous function
$f_K(w)=\min\{|w|^{-\kappa_\alpha},K\}$, with $f_K(0)=K$.
Then
\[
 \mathbb E\int f_K\,d\mu_{Y_n}
       \le C_\alpha+K\mathbb P(\mathcal G_n^c).
\]
Use Theorem~\ref{thm:interior-convergence}, first with fixed $K$, to
pass to the limit in the bounded integrals, and then let $K\to\infty$.
Monotone convergence proves \eqref{geo:eq-negative-radial}.
It follows that $\mu_\alpha^{\mathrm{SN}}(\{0\})=0$ and
$F_\alpha(r)\le C_\alpha r^{\kappa_\alpha}$.

Rotational invariance also rules out every nonzero point atom.
Indeed, for each $r>0$, if
$\mu_\alpha^{\mathrm{SN}}(\{|z|=r\})>0$, then the normalized
restriction of $\mu_\alpha^{\mathrm{SN}}$ to the circle
$\{|z|=r\}$ is a rotation-invariant probability measure.
By uniqueness of normalized Haar measure on the circle
(see, e.g., \cite[Theorem~2.20]{FollandAHA}), it is normalized
arc length and hence has no atoms. This proves Theorem  \ref{thm:limit-properties} (iii).
\end{proof}

\subsection{No mass on the unit circle}
\label{geo:subsec-boundary}

We first record a lower tail estimate for normalized stable averages.
%It will be applied to child resolvent variables, not to independent
normalizers substituted into the model.

\begin{lemma}[Thinning a normalized stable average]
\label{geo:lem-thinning}
Let $(\pi_k)$ have the law in Lemma~\ref{geo:lem-stable}, and let
$(V_k)$ be i.i.d. nonnegative random variables independent of these
weights. Suppose that $\mathbb P(V_1\ge a)\ge p_0>0$ for some $a>0$.
Then, for each $0<s<\beta$, there is
$C=C(\beta,s,a,p_0)<\infty$ such that
\begin{equation}\label{geo:eq-thinning-negative}
 \mathbb E\left(\sum_k\pi_kV_k\right)^{-s}\le C.
\end{equation}
\end{lemma}

\begin{proof}
Write
  $\pi_k=\frac{J_k}{S}$ and $S=\sum_kJ_k$,
where $(J_k)$ is the Poisson cloud from
Lemma~\ref{geo:lem-stable}. Let
\[
  I_k=\mathbf1_{\{V_k\ge a\}},
  \qquad
  p=\mathbb P(V_1\ge a)\ge p_0.
\]
Then
\[
  \sum_k\pi_kV_k
  \ge
  a\sum_k\pi_kI_k
  =
  a\,\frac{S_1}{S_0+S_1},
\]
where
  $S_1=\sum_kJ_kI_k$
and   $S_0=\sum_kJ_k(1-I_k)$.

By the Poisson thinning theorem
\cite[Corollary~5.9]{LastPenrose},
the retained and discarded point processes are independent Poisson
processes with intensities
\[
  p\beta x^{-1-\beta}\,dx
  \qquad \text{and } \qquad
  (1-p)\beta x^{-1-\beta}\,dx,
\]
respectively. Hence $S_1$ and $S_0$ are independent positive
$\beta$-stable random variables with Laplace transforms
  $\mathbb Ee^{-uS_1}
  =e^{-pK_\beta u^\beta}$ and $\mathbb Ee^{-uS_0}
  =e^{-(1-p)K_\beta u^\beta}$.

Therefore, for $0<s<\beta<1$,
\begin{align*}
  \mathbb E\left(\sum_k\pi_kV_k\right)^{-s}
  &\le
  a^{-s}\,
  \mathbb E\left(\frac{S_0+S_1}{S_1}\right)^s\\
  &=a^{-s}\,
  \mathbb E\left(1+\frac{S_0}{S_1}\right)^s\\
  &\le
  a^{-s}\left(
    1+\mathbb E S_0^s\,\mathbb E S_1^{-s}
  \right),
\end{align*}
where we used $(1+x)^s\le1+x^s$ for $0<s<1$ and the independence of
$S_0$ and $S_1$.
By the stable moment identities of
Lemma~\ref{geo:lem-stable},
\[
  \mathbb ES_0^s
  =
  \bigl((1-p)K_\beta\bigr)^{s/\beta}
  \frac{\Gamma(1-s/\beta)}{\Gamma(1-s)}
\quad \text{ and } \quad
  \mathbb ES_1^{-s}
  =
  \bigl(pK_\beta\bigr)^{-s/\beta}
  \frac{\Gamma(1+s/\beta)}{\Gamma(1+s)},
\]
which implies
\[
  \mathbb ES_0^s\,\mathbb ES_1^{-s}
  =
  \left(\frac{1-p}{p}\right)^{s/\beta}
  \frac{\Gamma(1-s/\beta)\Gamma(1+s/\beta)}
       {\Gamma(1-s)\Gamma(1+s)}.
\]
Since $p\ge p_0>0$, the right-hand side is bounded by a constant
depending only on $\beta,s,p_0$. Thus
\[
  \mathbb E\left(\sum_k\pi_kV_k\right)^{-s}
  \le C(\beta,s,a,p_0),
\]
which proves \eqref{geo:eq-thinning-negative}.
\end{proof}

\begin{proposition}[Theorem \ref{thm:limit-properties}(iv)]
\label{geo:prop-boundary}
Fix $r_1\in(\chi_\alpha^{1/4},1)$ and $0<s<\beta$.
There is $C_{\alpha,r_1,s}<\infty$ such that
\begin{equation}\label{geo:eq-boundary-tail}
 1-F_\alpha(r)
       \le C_{\alpha,r_1,s}(1-r^2)^{s/(1+s)},
                       \qquad r_1\le r<1.
\end{equation}
Consequently $\mu_\alpha^{\mathrm{SN}}\{|w|=1\}=0$.
\end{proposition}

\begin{proof}
Fix $r\in[r_1,1)$ and use $q=\mathbb EQ_o(r,t)$ as above.
By \eqref{geo:eq-t-over-q}, for sufficiently small $t$ (depending
on $r$),
\[
 \mathbb E\frac{H_+}{q}=1-\frac tq\ge\frac12.
\]
Furthermore, \eqref{geo:eq-relative-variance} and $H_+\le Q/r^2$ imply
\begin{equation}\label{geo:eq-H-relative-second}
 \mathbb E\left(\frac{H_+}{q}\right)^2
       \le\frac1{r^4-\chi_\alpha}
       \le\frac1{r_1^4-\chi_\alpha}.
\end{equation}
Put $\delta_*=r_1^4-\chi_\alpha>0$ and $V=H_+/q$.
Splitting at $1/4$ and applying Cauchy--Schwarz,
\[
 \frac12\le\mathbb EV
       \le\frac14+\delta_*^{-1/2}
                       \mathbb P(V\ge1/4)^{1/2}.
\]
Rearranging yields
\begin{equation}\label{geo:eq-good-child}
 \mathbb P(H_+/q\ge1/4)\ge\delta_*/16.
\end{equation}
The source recursion gives
\[
 \frac Qq=\frac tq+\sum_k\pi_k\frac{H_{+,k}}q
             \ge\sum_k\pi_k\frac{H_{+,k}}q.
\]
Apply Lemma~\ref{geo:lem-thinning}, using the independence from
Lemma~\ref{geo:lem-source-trace}, with threshold $1/4$ and success
probability at least $\delta_*/16$. This yields
\begin{equation}\label{geo:eq-Q-relative-negative}
 \mathbb E(Q/q)^{-s}\le C_{\alpha,r_1,s}.
\end{equation}
The constant does not depend on $r\in[r_1,1)$ or on sufficiently
small $t$ for that fixed $r$. The permitted upper bound on $t$ need
not be uniform as $r\uparrow1$; the limit in $t$ is taken first.

Let $A=Q/q$ and $D=D_{r,t}$. By \eqref{geo:eq-QD},
\[
 \mathbb E(AD)=1-r^2+r^2t/q.
\]
H\"older's inequality, $0\le D\le1$, and
\eqref{geo:eq-Q-relative-negative} give
\begin{align*}
 \mathbb ED
 &=\mathbb E\bigl[(AD)^{s/(1+s)}
                         (DA^{-s})^{1/(1+s)}\bigr]\\
 &\le(\mathbb EAD)^{s/(1+s)}
                         (\mathbb EDA^{-s})^{1/(1+s)}\\
 &\le C_{\alpha,r_1,s}
                         (1-r^2+r^2t/q)^{s/(1+s)}.
\end{align*}
At a continuity radius, let $t\downarrow0$ and use
\eqref{geo:eq-tail-recovery} and \eqref{geo:eq-t-over-q}.
This proves \eqref{geo:eq-boundary-tail}; right-continuity extends it
to the other radii. Since the measure is supported in the closed disk,
\[
 \mu_\alpha^{\mathrm{SN}}\{|w|=1\}
          =\lim_{r\uparrow1}(1-F_\alpha(r))=0,
\]
which proves Theorem \ref{thm:limit-properties} (iv).
\end{proof}

\subsection{The circular endpoint as \texorpdfstring{$\alpha\uparrow2$}{alpha increases to two}}
\label{geo:subsec-circular}

\begin{proposition}[Theorem \ref{thm:parameterendpoints} (i)]
\label{geo:prop-circular}
If $\chi_\alpha^{1/4}<r<1$, then
\begin{equation}\label{geo:eq-circular-pointwise}
 |F_\alpha(r)-r^2|
       \le\sqrt{\frac{\chi_\alpha}{r^4-\chi_\alpha}}.
\end{equation}
If $0<\chi_\alpha\le1/4$, then
\begin{equation}\label{geo:eq-circular-uniform}
 \sup_{r\ge0}|F_\alpha(r)-\min\{r^2,1\}|
       \le3\chi_\alpha^{1/4}.
\end{equation}
Moreover, $\mu_\alpha^{\mathrm{SN}}\Rightarrow\mu_{\mathrm{circ}}$
as $\alpha\uparrow2$.
\end{proposition}

\begin{proof}
Let $D=D_{r,t}$ and $q=\mathbb EQ_o(r,t)$. By
\eqref{geo:eq-QD}, we have
$ \mathbb E\left[\frac Qq(1-D)\right]
       =r^2(1-t/q)$.
Since $F_{\alpha,t}(r)=\mathbb E(1-D)$, Cauchy--Schwarz gives
\begin{align*}
 \left|F_{\alpha,t}(r)-r^2(1-t/q)\right|
 &=\left|\mathbb E\left[(1-Q/q)(1-D)\right]\right|\\
 &\le\left(\mathbb E|1-Q/q|^2\right)^{1/2}
 \le\sqrt{\frac{\chi_\alpha}{r^4-\chi_\alpha}}.
\end{align*}
The last bound is uniform in $t>0$.
Let $t\downarrow0$ at a continuity radius, using
\eqref{geo:eq-t-over-q} and Lemma~\ref{geo:lem-radial-recovery}.
This proves \eqref{geo:eq-circular-pointwise}. As before,
right-continuity extends it to every radius in the stated interval, see Proposition \ref{geo:prop-outer}.

Write $\chi=\chi_\alpha\le1/4$ and set $r_0=\chi^{1/8}$.
For $r_0\le r<1$,
\[
 |F_\alpha(r)-r^2|
 \le\sqrt{\frac{\chi}{\chi^{1/2}-\chi}}
 \le\sqrt2\,\chi^{1/4}.
\]
For $0\le r<r_0$, monotonicity yields
\[
 0\le F_\alpha(r)\le F_\alpha(r_0)
         \le(1+\sqrt2)\chi^{1/4},\qquad
 0\le r^2\le\chi^{1/4}.
\]
The absolute difference is therefore at most
$(1+\sqrt2)\chi^{1/4}<3\chi^{1/4}$. For $r\ge1$ both distribution
functions equal one, by the support inclusion. This proves
\eqref{geo:eq-circular-uniform}.

Finally, if $f$ is bounded and continuous on $\mathbb C$, let
\[
 \overline f(r)=\frac1{2\pi}\int_0^{2\pi}f(re^{i\vartheta})\,d\vartheta.
\]
Rotational invariance gives
$\int f\,d\mu_\alpha^{\mathrm{SN}}=\int_{[0,1]}\overline f\,dF_\alpha$.
The radial CDF convergence in \eqref{geo:eq-circular-uniform} gives
\[
 \int f\,d\mu_\alpha^{\mathrm{SN}}
       \longrightarrow\int_0^1\overline f(r)\,2r\,dr
       =\frac1\pi\int_{|z|\le1}f(z)\,d^2z.
\]
This is weak convergence to the circular law, which proves Theorem \ref{thm:parameterendpoints}(i). 
\end{proof}

\subsection{The endpoint at zero via a one-hot comparison}
\label{geo:subsec-zero}

Universality allows us to work throughout this subsection with the
symmetric Pareto representative \eqref{geo:eq-pareto-representative}
for each fixed $\alpha$. It is convenient to use the original
row-normalized matrices $Y_n$ rather than their transposes.
For the $i$-th row of $Y_n$, let $I_{i,n}$ be the column index $j\in [n]$ that maximizes $|Y_{ij}|$ and put
\begin{equation}\label{geo:eq-onehot}
 (B_n)_{ij}=\operatorname{sgn}(Y_{iI_{i,n}})
                              \mathbf1_{\{j=I_{i,n}\}}.
\end{equation}
There are almost surely no ties for this representative. The indices
$I_{i,n}, i\in [n]$ are independent and uniform on $[n]$, and the selected
signs are independent Rademacher variables, independent of the indices.
In \cite{franccois2026characteristic}, the matrix $B_n$ was called the signed one-hot matrix associated to $Y_n$. 

\begin{lemma}[One-hot approximation]
\label{geo:lem-onehot-error}
Let $\varepsilon_{n,\alpha} :=n^{-1}\|Y_n-B_n\|_F^2$. The signed one-hot matrix $B_n$ satisfies
\begin{equation}\label{geo:eq-onehot-error}
 \limsup_{n\to\infty}\mathbb E \varepsilon_{n,\alpha}\le\alpha.
\end{equation}
\end{lemma}

\begin{proof}
As in the proof of \cite[Lemma 2.10]{franccois2026characteristic}, one can compute that
\[
 \varepsilon_{n,\alpha}
  =\frac2n\sum_{i=1}^n(1-\max_j|Y_{ij}|).
\]
For a unit vector $y$, writing $m=\max_j|y_j|\le1$ gives
$\sum_j|y_j|^4\le m^2\le m$. Hence
\begin{equation}\label{geo:eq-error-fourth}
 \mathbb E\varepsilon_{n,\alpha}
       \le2\Big(1-\mathbb E\sum_{j=1}^n|Y_{1j}|^4\Big).
\end{equation}
From \cite[Lemma 4.1(a)]{franccois2026characteristic}, we know that 
$\lim_{\nto} \mathbb E\sum_{j=1}^n|Y_{1j}|^4 = 1-\alpha/2$,
which in conjunction with \eqref{geo:eq-error-fourth} concludes the proof.
\end{proof}

%f_n^k(i):= \arg \max_{j\in [n]} \big|(B_n^k)_{ij}\big|
\begin{lemma}[Powers of one-hot matrices]
\label{geo:lem-mapping}
(a) For positive integers $k$, let 
$$f_n(i):= \arg \max_{j\in [n]} \big|(B_n)_{ij}\big| =I_{n,i} \qquad  \text{ and } \qquad f_n^{\circ k}:=f_n\circ \cdots \circ f_n$$
be the $k$-th iterate of $f_n$. Then 
\begin{equation}\label{geo:eq-onehot-rank}
 \|B_n^k\|_F^2=n \qquad \text{ and } \qquad
 \operatorname{rank}(B_n^k)=|f_n^{\circ k}([n])|\,, \quad k\ge1.
\end{equation}
(b) Define $q_0=1$ and
 $q_{k+1}=1-e^{-q_k}$, $k\ge 1$.
Then it holds for $k\ge 1$ that $q_k\le\frac2{k+2}$,
\begin{equation}\label{geo:eq-rank-limit}
 \lim_{n\to\infty}
 \mathbb E\frac{\operatorname{rank}(B_n^k)}n=q_k \qquad \text{ and } 
 \qquad
 \frac{\operatorname{rank}(B_n^k)}n \longrightarrow  q_k \quad \text{ a.s.}
\end{equation}
\end{lemma}

\begin{proof}
The $i$-th row of $B_n$ has a single nonzero entry of modulus one, in
column $f_n(i)$. One sees that the $i$-th row of $B_n^2$ also has only one nonzero entry of modulus one located in column $f_n(f_n(i))=f_n^{\circ 2}(i)$. 
Iterating shows that the $i$-th row of $B_n^k$
again has a single nonzero entry of modulus one, now in column
$f_n^{\circ k}(i)$. Hence
 $\|B_n^k\|_F^2=n$,
and the column space of $B_n^k$ is spanned by the columns indexed by
$f_n^{\circ k}([n])$. Since those columns are linearly independent, we obtain 
\[
 \operatorname{rank}(B_n^k)=|f_n^{\circ k}([n])|.
\]

For the uniform random mapping $f_n:[n]\to[n]$, the classical
random-mapping estimate of
\cite[Theorem~2]{FlajoletOdlyzko}
gives, for every fixed $k$,
\[
 \frac1n\mathbb E|f_n^{\circ k}([n])|
 \longrightarrow 1-\tau_k, \qquad \nto\,
\]
where $\tau_0=0$ and $\tau_{k+1}=e^{-1+\tau_k}$.
Setting $q_k=1-\tau_k$ gives
\[
 q_0=1,\qquad
 q_{k+1}=1-e^{-q_k},
\]
and proves the first assertion in \eqref{geo:eq-rank-limit}.

For the almost surely convergence since for every fixed $k$ changing one image of the random mapping changes \(|f_n^{\circ k}([n])|\) by at most \(2k\). So bounded differences inequality yields
    \[
\mathbb P\!\left(
\left|
\frac{|f_n^{\circ k}([n])|}{n}
-
\mathbb E\frac{|f_n^{\circ k}([n])|}{n}
\right|>\delta
\right)
\le
2\exp\!\left(-\frac{\delta^2n}{2k^2}\right),
\]
which is summable in $n$. Therefore, the second claim in \eqref{geo:eq-rank-limit} follows by Borel–Cantelli.

Finally, for $0\le x\le1$, we have
 $1-e^{-x}\le\frac{x}{1+x/2}$.
Since $0<q_k\le1$, this implies
\[
 \frac1{q_{k+1}}
 \ge\frac1{q_k}+\frac12.
\]
Starting from $q_0=1$ and iterating, we get
\[
 \frac1{q_k}\ge1+\frac k2,
\]
which concludes the proof.
\end{proof}

For completeness, the next lemma provides the Schatten inequalities
used in the power comparison, including the exponents below one.
For an $n\times n$ matrix $M$ and $p>0$, we introduce the normalized Schatten $p$-norm 
\[
 \|M\|_{p,n}=\left(\frac1n\sum_j s_j(M)^p\right)^{1/p}
\]
which is a matrix norm for every $p\in [1,\infty]$. For $p\in (0,1)$, it is a quasi-norm.

\begin{lemma}[Schatten inequalities]
\label{geo:lem-schatten}
For $0<p\le1$ and square matrices $A,B$,
\begin{equation}\label{geo:eq-p-subadditive}
 \|A+B\|_{p,n}^p
 \le\|A\|_{p,n}^p+\|B\|_{p,n}^p.
\end{equation}
For square matrices $A_1,\ldots,A_k$,
\begin{equation}\label{geo:eq-holder}
 \|A_1\cdots A_k\|_{2/k,n}
 \le\prod_{\ell=1}^k\|A_\ell\|_{2,n}.
\end{equation}
For every $p>0$,
\begin{equation}\label{geo:eq-weyl-majorant}
 \sum_j|\lambda_j(A)|^p
 \le\sum_j s_j(A)^p.
\end{equation}
\end{lemma}
\begin{proof}
Equation~\eqref{geo:eq-p-subadditive} is the
Rotfel'd--Thompson inequality applied to the increasing concave
function $x\mapsto x^p$; see \cite[Theorem~3]{Thompson76}.
Equation~\eqref{geo:eq-holder} is the Schatten  H\"older inequality
iterated with all factor exponents equal to $2$, and
\eqref{geo:eq-weyl-majorant} is Weyl's majorant theorem; see
\cite[Chapter~3]{HornJohnsonTopics}.
\end{proof}

\begin{lemma}[Rank--power comparison]
\label{geo:lem-rank-power}
Suppose $A,B$ are $n\times n$ matrices with
$\|A\|_F^2=\|B\|_F^2=n$, and $B$ is a signed one-hot matrix.
Then, for every integer $k\ge2$,
\begin{equation}\label{geo:eq-rank-power}
 \frac1n\sum_j|\lambda_j(A)|^2
 \le\left(\frac{\operatorname{rank}(B^k)}n\right)^{1-1/k}
       +k\left(\frac{\|A-B\|_F^2}{n}\right)^{1/k}.
\end{equation}
\end{lemma}

\begin{proof}
Set $p=2/k\le1$. By \eqref{geo:eq-weyl-majorant},
\[
 \frac1n\sum_j|\lambda_j(A)|^2
   =\frac1n\sum_j|\lambda_j(A^k)|^{2/k}
   \le\|A^k\|_{p,n}^p.
\]
Recall the algebraic identity
 $A^k=B^k+\sum_{\ell=0}^{k-1}
                     A^\ell(A-B)B^{k-1-\ell}$.
Each summand is a product of exactly $k$ factors. Since
$\|A\|_{2,n}=\|B\|_{2,n}=1$, \eqref{geo:eq-holder} gives
\[
 \|A^\ell(A-B)B^{k-1-\ell}\|_{2/k,n}^{2/k}
       \le\left(\frac{\|A-B\|_F^2}{n}\right)^{1/k}.
\]
Let $d=\operatorname{rank}(B^k)$. The matrix $B^k$ is again signed
one-hot, so $\sum_j s_j(B^k)^2=n$. Scalar H\"older inequality gives
\[
 \|B^k\|_{2/k,n}^{2/k}
  =\frac1n\sum_{j=1}^d s_j(B^k)^{2/k}
  \le\frac1n d^{1-1/k}
                   \left(\sum_j s_j(B^k)^2\right)^{1/k}
  =\left(\frac dn\right)^{1-1/k}.
\]
By \eqref{geo:eq-p-subadditive}, these estimates prove
\eqref{geo:eq-rank-power}.
\end{proof}

\begin{proposition}[Key to Theorem \ref{thm:parameterendpoints}(ii)]
\label{geo:prop-zero-endpoint}
For every integer $k\ge2$,
\begin{equation}\label{geo:eq-second-moment-zero}
 \int_{\mathbb C}|w|^2\,\mu_\alpha^{\mathrm{SN}}(dw)
   \le\left(\frac2{k+2}\right)^{1-1/k}+k\alpha^{1/k}.
\end{equation}
\end{proposition}

\begin{proof}
Keep $\alpha>0$ and $k$ fixed and let $n\to\infty$.
Lemmas~\ref{geo:lem-onehot-error} and~\ref{geo:lem-mapping} give
\[
 \limsup_{n\to\infty}
 \mathbb E\frac1n\sum_{j=1}^n|\lambda_j(Y_n)|^2
 \le q_k^{1-1/k}+k\alpha^{1/k}.
\]
Since $\mu_\alpha^{\mathrm{SN}}$ is supported in
$\overline{\mathbb D}$, the bounded continuous function
\[
 f(w)=\min\{|w|^2,1\}
\]
satisfies
\[
 \int_{\mathbb C}f(w)\,\mu_\alpha^{\mathrm{SN}}(dw)
 =
 \int_{\mathbb C}|w|^2\,\mu_\alpha^{\mathrm{SN}}(dw).
\]
By Theorem~\ref{thm:interior-convergence},
\[
 \int f\,d\mu_{Y_n}
 \longrightarrow
 \int |w|^2\,\mu_\alpha^{\mathrm{SN}}(dw)
\]
in probability, and hence also in expectation, since
$0\le\int f\,d\mu_{Y_n}\le1$. Moreover,
\[
 \int f\,d\mu_{Y_n}
 \le \frac1n\sum_{j=1}^n|\lambda_j(Y_n)|^2.
\]
Therefore,
\[
 \int_{\mathbb C}|w|^2\,\mu_\alpha^{\mathrm{SN}}(dw)
 \le q_k^{1-1/k}+k\alpha^{1/k}\le\left(\frac2{k+2}\right)^{1-1/k}+k\alpha^{1/k},
\]
where Lemma~\ref{geo:lem-radial-recovery}(b) was used for the last step. This proves \eqref{geo:eq-second-moment-zero}.
\end{proof}

\begin{proof}[Completion of the proof of Theorem \ref{thm:parameterendpoints}(ii)]
For fixed $k$, let $\alpha\downarrow0$ in
\eqref{geo:eq-second-moment-zero}. We get
\[\limsup_{\alpha\downarrow 0} \int_{\mathbb C}|w|^2\,\mu_\alpha^{\mathrm{SN}}(dw)
 \le  \left(\frac2{k+2}\right)^{1-1/k}.\]
Now let $k\to\infty$ to obtain
\[
 \int_{\mathbb C} |w|^2\,\mu_\alpha^{\mathrm{SN}}(dw)\longrightarrow0\,, \qquad \alpha\downarrow 0.
\]
For every $\varepsilon>0$, we conclude that
\[
 \mu_\alpha^{\mathrm{SN}}\{|w|>\varepsilon\}
 \le\varepsilon^{-2}
           \int_{\mathbb C} |w|^2\,\mu_\alpha^{\mathrm{SN}}(dw)
 \longrightarrow0 \,, \qquad \alpha\downarrow 0,
\]
which proves $\mu_\alpha^{\mathrm{SN}}\Rightarrow\delta_0$ as
$\alpha\downarrow0$.
\end{proof}

\section{Proofs at the exact endpoints}\label{sec:endpoint-proofs}
 We now prove the two exact endpoint
theorems stated in Section~\ref{sec:model-main-results}.  Their proofs use
different mechanisms.  At $\alpha=2$ the row norm is relatively stable and
self-normalization is asymptotically deterministic; at $\alpha=0$ one
coordinate asymptotically carries the whole row mass and the matrix behaves
asymptotically like a signed one-hot matrix.

\subsection{The normal-domain endpoint: \texorpdfstring{$\alpha=2$}{alpha = 2}}\label{end:subsec-two}

In this subsection we will prove Theorem \ref{thm:exact-two} under the assumptions therein. We start by choosing a deterministic sequence $c_n\to\infty$ such that
\begin{equation}\label{eq:model-normal-domain-scale}
  \frac{n\ell_2(c_n)}{c_n^2}\longrightarrow1.
\end{equation}
Such a sequence exists by asymptotic inversion, since
$x\mapsto x^2/\ell_2(x)$ is regularly varying with index two;
see \cite[Theorem~1.5.12]{BGT87}.
%The scale $c_n$ is the concentration scale of the Euclidean row norm, see Lemma \ref{end:lem-row-two}.
When $\mathbb E|X|^2=\sigma^2\in(0,\infty)$, one may take
$c_n=\sigma\sqrt n$. %In general, $c_n$ need not be asymptotic to the
extreme-value quantile associated with the tail of $|X|$.

\begin{lemma}[Relative stability of the row norm]
\label{end:lem-row-two}
Let $L_n=\left(\sum_{j=1}^n|X_j|^2\right)^{1/2}$,
where $X_1,\ldots,X_n$ are independent copies of $X$. Then, as $\nto$,
\begin{equation}\label{end:eq-row-relative-two}
  \frac{L_n}{c_n}\longrightarrow1
  \quad\text{in probability},
  \qquad
  \frac{\max_{j\le n}|X_j|}{c_n}\longrightarrow0
  \quad\text{in probability}.
\end{equation}
Moreover, for every $0<p<2$,
\begin{equation}\label{end:eq-row-Lp-two}
  \sup_{n\ge n_0}
  \mathbb E\left(\frac{L_n}{c_n}\right)^p<\infty,
  \qquad
  \mathbb E\left|
    \frac{L_n}{c_n}-1
  \right|^p\longrightarrow0.
\end{equation}
\end{lemma}

\begin{proof}
Apply the relative-stability theorem for sums of nonnegative random
variables with slowly varying truncated first moment to
$Z_j=|X_j|^2$ and $b_n=c_n^2$; see
\cite[Chapter~8]{BGT87}. This gives
  $L_n^2/c_n^2\longrightarrow1$ in probability.
The normal-domain characterization
\[
  x^2\mathbb P(|X|>x)=o(\ell_2(x))
\]
and a union bound give
$\max_{j\le n}|X_j|/c_n\to0$ in probability.

Moreover,
\[
  \mathbb P(L_n>u)
  \le n\mathbb P(|X|>u)+\frac{n\ell_2(u)}{u^2}.
\]
Taking $u=c_nT$ and applying Potter's bound yields, for every
$\varepsilon>0$,
\[
  \mathbb P(L_n>c_nT)\le C_\varepsilon T^{-2+\varepsilon},
  \qquad T\ge1.
\]
Tail integration gives the uniform $p$-th moment bound for $p<2$.
Together with convergence in probability, this proves the remaining
$L^p$ assertions by uniform integrability.
\end{proof}

\begin{proposition}[Comparison with deterministic normalization]
\label{end:prop-comparison-two}
For every $z\in\C$,
\begin{equation}\label{end:eq-mirsky-two}
 \frac1n\sum_{j=1}^n
 \left|s_j(Y_n-zI)-s_j(A_n-zI)\right|
 \le\frac1n\sum_{i=1}^n\left|\frac{D_{i,n}}{c_n}-1\right|.
\end{equation}
where $A_n=c_n^{-1}X_n$ and $D_{i,n}$ is as in \eqref{eq:model-spectral-measures}. Consequently, as $\nto$,
\begin{equation}\label{end:eq-BL-two}
 \Ee\sup_{z\in\C}
 d_{\mathrm{BL}}\bigl(\nu_{Y_n-zI},\nu_{A_n-zI}\bigr)
 \longrightarrow0.
\end{equation} 
\end{proposition}

\begin{proof}
By Mirsky's inequality,
\[
 \sum_{j=1}^n
 |s_j(Y_n-zI)-s_j(A_n-zI)|
 \le \|Y_n-A_n\|_{S_1},
 \qquad
 \text{ where } \|M\|_{S_1}:=\sum_j s_j(M).
\]
Writing $M_{i\bullet}$ for the $i$-th row and using the triangle
inequality together with
$\|uv^*\|_{S_1}=\|u\|_2\|v\|_2$, yields
\[
 \|Y_n-A_n\|_{S_1}
 \le\sum_{i=1}^n\|(Y_n-A_n)_{i\bullet}\|_2
 =\sum_{i=1}^n\left|1-\frac{D_{i,n}}{c_n}\right|.
\]
See \cite[Chapter~3]{HornJohnsonTopics}. Divide by $n$, take expectations, and use
Lemma~\ref{end:lem-row-two}. The bounded-Lipschitz assertion follows
by testing the paired singular values against $1$-Lipschitz functions.
\end{proof}

The next lemma identifies the reference law.  We state the Gaussian
shifted singular value limit as the standard input from the finite variance
circular law theory; see, for instance, \cite[Corollary~4.10]{bordenave2012around}.

\begin{lemma}[Shifted reference singular value limit]
\label{end:lem-reference-two}
For every fixed $z\in\mathbb C$,
\begin{equation}\label{end:eq-reference-two}
  \nu_{A_n-zI}\Rightarrow\nu_z^{\mathrm G}
  \qquad\text{in probability}.
\end{equation}
\end{lemma}

\begin{proof}
The claim follows by Corollary 2.10 of \cite{AdamczakDS}.
\end{proof}

\begin{proof}[Proof of \eqref{eq:exact-two-singular-values}, the shifted singular value part of
Theorem~\ref{thm:exact-two}]
Follows by Proposition~\ref{end:prop-comparison-two} and 
Lemma~\ref{end:lem-reference-two}.
\end{proof}

To pass from shifted singular values to eigenvalues we only need a
normal-domain replacement for the high-codimension inverse-distance estimate.
The deterministic negative-moment inequality, the row-to-hyperplane
argument, and the Hermitization step have already been proved in
Sections~\ref{ex:subsec-invertibility}--\ref{ex:subsec-log-ui}.

\begin{lemma}[Normal-domain high-codimension inverse distances]
\label{end:lem-distance-two}
Let $V\subset\C^n$ be a deterministic complex linear subspace of
codimension $d\ge n^{7/8}$, and set
  $D_V=\dist((X_1,\ldots,X_n)^\top,V)$.
Then there exists an event $E_V$ such that
  $1- \Pp(E_V)\le Ce^{-cn^{1/8}}$,
and, for every $a>0$,
\begin{equation}\label{end:eq-distance-raw-two}
  \Ee[D_V^{-a};E_V]
  \le C_a\left(\frac{n}{c_nd}\right)^a.
\end{equation}
Moreover, if $W$ is independent of $(X_1,\ldots,X_n)$,
  $V=W+\Span\{e_i\}$
and $U=(X_1,\ldots,X_n)^\top/\|(X_1,\ldots,X_n)^\top\|_2$, then for every $0<p<1$,
\begin{equation}\label{end:eq-distance-normalized-two}
  \Ee[\dist(U-ze_i,W)^{-p};E_V\mid W]
  \le C_p\left(\frac nd\right)^p,
\end{equation}
uniformly over deterministic $z\in\C$.
\end{lemma}

\begin{proof}
We follow the proof of the distance-to-subspace estimate
\cite[Lemma~4.13]{bordenave2012around}, retaining only the modifications caused by
the normal-domain scale.

For $u_0\le u\le n^{1/8}$, let
  $b:=c_n\frac dn\,u^{-1}$ and
condition on the retained set
\[
  I=\{j:|X_j|\le b\}.
\]
After adjoining the constant vector to the projected comparison
subspace, the retained coordinates may be centered exactly as in the
proof of Lemma~4.13 of \cite{bordenave2012around}. Their conditional
variance is
\[
  \tau_b^2
  :=
  \Ee\bigl[
    |X-\Ee(X\mid |X|\le b)|^2
    \,\bigm|\,
    |X|\le b
  \bigr]
  \sim\ell_2(b).
\]

The normal-domain estimates and Potter's bound give, uniformly in the
stated range,
\[
  \Ee|I^c|
  =n\Pp(|X|>b)
  \le C\left(\frac nd\right)^2u^2=o(d),
\]
and
\[
  \frac{d\tau_b^2}{b^2}\ge cu,
  \qquad
  \sqrt{d\tau_b^2}
  \ge c\,c_n\frac dn\,u^{-1/2}.
\]
A Chernoff bound therefore shows that
$\Pp(|I^c|>d/4)\le e^{-cd}$.
The median-concentration argument in the proof of
\cite[Lemma~4.13]{bordenave2012around}, applied conditionally on $I$,
then yields
\begin{equation}\label{end:eq-distance-smallball-two}
  \Pp\left(
    D_V\le c_*c_n\frac dn\,u^{-1/2}
  \right)
  \le e^{-cd}+Ce^{-cu}.
\end{equation}

Taking $u=n^{1/8}$ defines $E_V$ and proves
$1- \Pp(E_V)\le Ce^{-cn^{1/8}}$. Integrating
\eqref{end:eq-distance-smallball-two} after the substitution $u=v^2$
gives \eqref{end:eq-distance-raw-two} for every fixed $a>0$.

Finally,
\[
  \dist(U-ze_i,W)
  \ge\dist(U,V)
  =\frac{D_V}{\|(X_1,\ldots,X_n)^\top}\|_2.
\]
For $0<p<1$, Cauchy--Schwarz,
Lemma~\ref{end:lem-row-two} at exponent $2p<2$, and
\eqref{end:eq-distance-raw-two} at exponent $2p$ give
\[
\begin{aligned}
  \Ee[\dist(U-ze_i,W)^{-p};E_V\mid W]
  &\le
  \bigl(\Ee\|X\|_2^{2p}\bigr)^{1/2}
  \bigl(\Ee[D_V^{-2p};E_V\mid W]\bigr)^{1/2}\\
  &\le C_p\left(\frac nd\right)^p.
\end{aligned}
\]
\end{proof}
\begin{proposition}[Logarithmic uniform integrability at the normal endpoint]
\label{end:prop-ui-two}
Assume in addition that $X$ has a bounded density.  For every fixed
$z\in\C$, $\log$ is uniformly integrable in probability with respect to
$(\nu_{Y_n-zI})_n$.
\end{proposition}

\begin{proof}
The polynomial least-singular value estimate
\eqref{ex:eq-invertibility} applies verbatim.  

For intermediate singular values repeat the row-deletion argument of
Section~\ref{ex:subsec-log-ui}, replacing Proposition~\ref{ex:prop-normalized-distance}
by Lemma~\ref{end:lem-distance-two}.  If
$i\ge 4n^{7/8}$ and $0<p<1$, the same deterministic fractional
negative-moment inequality gives, on a common event of probability tending
to one,
\begin{equation}\label{end:eq-intermediate-two}
 \Ee\bigl[s_{n-i}(Y_n-zI)^{-p}; E_V\bigr]
 \le C\left(\frac ni\right)^{1+p}.
\end{equation}
Choose $0<q<p$ with
$q(1+p)/p<1$, and also small enough that the bottom $O(n^{7/8})$
singular values are negligible under the polynomial least-singular value
bound.  Exactly as in the proof of Proposition~\ref{ex:prop-ui}, one obtains
\[
 \sup_n\Ee\left[
   \mathbf1_{G_n(z)}\int_0^\infty s^{-q}\nu_{Y_n-zI}(ds)
 \right]<\infty
\]
for events $G_n(z)$ with probability tending to one.  This controls the
negative logarithmic tail.  The positive logarithmic tail is still bounded
deterministically by \eqref{eq:model-unit-rows}.  Hence $\log$ is uniformly
integrable in probability.
\end{proof}

\begin{proof}[Completion of the proof of Theorem~\ref{thm:exact-two}]
The shifted singular value convergence was proved above.  Under bounded
density, Proposition~\ref{end:prop-ui-two} allows the same Hermitization
argument as in the proof of Theorem ~\ref{thm:interior-convergence}.  The limiting shifted
singular value family is the Gaussian family $\nu_z^{\mathrm G}$, whose
logarithmic potential is the circular-law potential.  Therefore,
  $\mu_{Y_n}\Rightarrow\mu_{\mathrm{circ}}$ in probability.
%At $z=0$, $\nu_0^{\mathrm G}$ is the quarter-circle law, giving the stated zero-shift singular value conclusion.
\end{proof}

\subsection{The slowly varying endpoint: \texorpdfstring{$\alpha=0$}{alpha = 0}}\label{end:subsec-zero}

We now assume that $X_{1,1}$ is real-valued and that 
$\Pp(|X_{1,1}|>t)$ is slowly varying.  Let $B_n$ be the signed one-hot matrix defined in \eqref{geo:eq-onehot}.  Recall the random-mapping and deterministic
rank--power estimates of Lemmas~\ref{geo:lem-mapping} and
\ref{geo:lem-rank-power}.

\begin{lemma}\cite[Theorem 3.2]{Darling52}
\label{end:lem-maximum-zero}
Let $Z_1,\ldots,Z_n$ be i.i.d.\ nonnegative random variables, whose
tail $\overline F(x)=\mathbb P(Z_1>x)$ is slowly varying at infinity. Put
  $S_n=\sum_{j=1}^n Z_j$ and  $M_n=\max_{1\le j\le n}Z_j$.
Then, as $\nto$, we have $M_n/S_n\longrightarrow1$ in probability and 
\begin{equation}\label{end:eq-maximum-zero}
  \mathbb E\frac{S_n-M_n}{M_n}\longrightarrow0.
\end{equation}
where the ratio is defined to be zero on $\{M_n=0\}$.
\end{lemma}

The one-hot approximation in Lemma~\ref{geo:lem-onehot-error} remains valid for $\alpha=0$. The next result is a stronger version of that statement for $\alpha=0$.
\begin{proposition}[One-hot approximation for $\alpha=0$]
\label{end:prop-frobenius-zero}
It holds
\begin{equation}\label{end:eq-frobenius-zero}
  \Delta_n^2:=\frac1n\|Y_n-B_n\|_F^2\longrightarrow0
  \qquad\text{almost surely and in }L^1.
\end{equation}
Moreover, uniformly in $z\in\C$,
\begin{equation}\label{end:eq-shift-zero}
  W_2\bigl(\nu_{Y_n-zI},\nu_{B_n-zI}\bigr)^2\le\Delta_n^2.
\end{equation}
\end{proposition}

\begin{proof}
We apply Lemma~\ref{end:lem-maximum-zero} with $Z_j=X_{ij}^2$.  On a nonzero
$i$-th row, we have
 $\max_j|Y_{ij}|=\sqrt{M_n/S_n}$,
and hence
\[
  \varepsilon_n:=2\Ee\bigl[1-\max_j|Y_{ij}|\bigr]\longrightarrow0.
\]
For a unit row $y$ and a maximizing coordinate $I$,
\[
 \|y-\operatorname{sgn}(y_I)e_I^{\mathsf T}\|_2^2
 =2(1-|y_I|).
\]
Therefore
\[
 \Delta_n^2=\frac2n\sum_{i=1}^n(1-\max_j|Y_{ij}|),
 \qquad \Ee\Delta_n^2=\varepsilon_n.
\]
For fixed $n$ the summands are independent and lie in $[0,2]$; Hoeffding's
inequality gives
\[
 \Pp(|\Delta_n^2-\varepsilon_n|>u)\le2e^{-nu^2/2}.
\]
Borel--Cantelli gives the almost-sure convergence, while boundedness gives
$L^1$ convergence.  Finally, the singular value Hoffman--Wielandt
inequality applied to $Y_n-zI$ and $B_n-zI$ proves
\eqref{end:eq-shift-zero}.
\end{proof}

\begin{proposition}[Eigenvalue collapse at the exact zero endpoint]
\label{end:prop-eigen-zero}
Under the hypotheses of Theorem~\ref{thm:exact-zero},
\begin{equation}\label{end:eq-eigen-zero}
 \frac1n\sum_{j=1}^n|\lambda_j(Y_n)|^2\longrightarrow0
 \qquad\text{almost surely and in }L^1.
\end{equation}
In particular, $\mu_{Y_n}\Rightarrow\delta_0$ almost surely.
\end{proposition}

\begin{proof}
Fix $k\ge2$.  Lemma~\ref{geo:lem-rank-power}, applied with
$A=Y_n$ and $B=B_n$, gives
\[
 \frac1n\sum_j|\lambda_j(Y_n)|^2
 \le\left(\frac{\rank(B_n^k)}n\right)^{1-1/k}
    +k(\Delta_n^2)^{1/k}.
\]
By Proposition~\ref{end:prop-frobenius-zero} and
Lemma~\ref{geo:lem-mapping}, for every fixed $k$,
\[
 \limsup_{n\to\infty}\frac1n\sum_j|\lambda_j(Y_n)|^2
 \le q_k^{1-1/k}\qquad\text{almost surely}.
\]
Intersect over the countable set of integers $k\ge2$ and then let
$k\to\infty$; since $q_k\le2/(k+2)$, the right-hand side tends to zero.
The left-hand side is bounded by $n^{-1}\|Y_n\|_F^2=1$, so convergence
also holds in $L^1$.  Markov's inequality for the empirical second moment
implies weak convergence to $\delta_0$.
\end{proof}

\begin{proposition}[Poisson squared-singular value limit]
\label{end:prop-poisson-zero}
Let
\[
 \Pi_1=e^{-1}\sum_{\ell=0}^\infty\frac1{\ell!}\delta_\ell,
 \qquad
 \nu_0^{\mathrm{SN}}
   =e^{-1}\sum_{\ell=0}^\infty\frac1{\ell!}\delta_{\sqrt\ell}.
\]
Then
\begin{equation}\label{end:eq-poisson-zero}
 H_{Y_n}\Rightarrow\Pi_1,
 \qquad
 \nu_{Y_n}\Rightarrow\nu_0^{\mathrm{SN}}
 \qquad\text{almost surely}.
\end{equation}
\end{proposition}

\begin{proof}
Let $f_n(i)$ be the column selected by row $i$ of $B_n$, and set
$N_{j,n}=\#\{i:f_n(i)=j\}$.  Since each row of $B_n$ has one entry of
modulus one,
\[
  B_n^*B_n=\diag(N_{1,n},\ldots,N_{n,n}),
  \qquad
  H_{B_n}=\frac1n\sum_{j=1}^n\delta_{N_{j,n}}.
\]
For every fixed $\ell\ge0$,
\[
 \Ee\frac1n\#\{j:N_{j,n}=\ell\}
 =\Pp(\operatorname{Binomial}(n,1/n)=\ell)
 \longrightarrow e^{-1}/\ell!.
\]
Changing one image $f_n(i)$ changes at most two occupancy counts, so
bounded differences and Borel--Cantelli give almost-sure convergence of
each empirical point mass.  Since the limiting masses sum to one,
$H_{B_n}\Rightarrow\Pi_1$ almost surely.  Both measures have first moment
one, and the same occupancy argument gives convergence in $W_1$.

To transfer the limit to $Y_n$, pair ordered singular values and use
Cauchy--Schwarz together with \eqref{end:eq-shift-zero} at $z=0$:
\[
 \begin{aligned}
 W_1(H_{Y_n},H_{B_n})
 &\le\frac1n\sum_j|s_j(Y_n)^2-s_j(B_n)^2|\\
 &\le \Delta_n
 \left(\frac1n\sum_j(s_j(Y_n)+s_j(B_n))^2\right)^{1/2}
 \le2\Delta_n.
 \end{aligned}
\]
Proposition~\ref{end:prop-frobenius-zero} therefore gives
$H_{Y_n}\Rightarrow\Pi_1$.  Finally,
$|\sqrt x-\sqrt y|^2\le|x-y|$ transfers the convergence through the
square-root map, proving the singular value statement.
\end{proof}

\begin{proof}[Completion of the proof of Theorem~\ref{thm:exact-zero}]
Proposition~\ref{end:prop-eigen-zero} proves the eigenvalue assertions,
and Proposition~\ref{end:prop-poisson-zero} proves the squared-singular-
value and singular value assertions.  Proposition~\ref{end:prop-frobenius-zero}
records the one-hot equivalence used in both parts. 
\end{proof}

\paragraph{Acknowledgments:}
This research was supported by the Swedish Research Council through VR-2023-03577 ``High-dimensional extremes and random matrix structures'' and by the Verg-Foundation. The authors thank Xuechun Hu for providing figure 1 and Quentin François for fruitful discussions.

\paragraph{AI Statement:}
The authors acknowledge the use of OpenAI's ChatGPT (GPT-5.6 Sol) during the preparation of this manuscript. The model was used as a research and writing aid, including for improving exposition, checking intermediate calculations, and suggesting relevant references. All mathematical arguments, statements, citations, and conclusions were written, reviewed and verified by the authors. The authors take full responsibility for the content of the manuscript, including any errors or omissions.

\bibliographystyle{alpha}

\bibliography{references}
\end{document}